\documentclass[11pt,article]{amsart}
\usepackage{graphicx,epsfig}
\usepackage{url}
\usepackage[english]{babel}
\usepackage[table]{xcolor}
\definecolor{cyan(process)}{rgb}{0.0, 0.72, 0.92}
\usepackage[ansinew]{inputenc}
\usepackage{a4wide}
\usepackage{graphicx}
\definecolor{columbiablue}{rgb}{0.61, 0.87, 1.0}
\usepackage{amsmath}
\usepackage{mathptmx}       
\usepackage{amssymb,amsthm,amsfonts,amstext,bbm}
\usepackage{amsmath}
\usepackage{enumerate}
\usepackage{a4wide}
\usepackage[active]{srcltx}
\usepackage{multicol}
\usepackage[charter]{mathdesign}
\definecolor{myblue}{rgb}{0.61, 0.87, 1.0}

\usepackage{amsmath}
\usepackage{empheq}

\newlength\mytemplen
\newsavebox\mytempbox

\makeatletter
\newcommand\mybluebox{%
    \@ifnextchar[
       {\@mybluebox}%
       {\@mybluebox[0pt]}}

\def\@mybluebox[#1]{%
    \@ifnextchar[
       {\@@mybluebox[#1]}%
       {\@@mybluebox[#1][0pt]}}

\def\@@mybluebox[#1][#2]#3{
    \sbox\mytempbox{#3}%
    \mytemplen\ht\mytempbox
    \advance\mytemplen #1\relax
    \ht\mytempbox\mytemplen
    \mytemplen\dp\mytempbox
    \advance\mytemplen #2\relax
    \dp\mytempbox\mytemplen
    \colorbox{myblue}{\hspace{1em}\usebox{\mytempbox}\hspace{1em}}}

\usepackage{lipsum}
\usepackage{eso-pic}
\usepackage{graphicx}

\definecolor{verdescuro}{RGB}{0,120,0}

\newcommand{\ola}[1]{\overleftarrow{#1}}
\newcommand{\ora}[1]{\overrightarrow{#1}}

\newtheorem{prop}{Proposition}[section]
\newtheorem{thm}[prop]{Theorem}

\newtheorem{definition}[prop]{Definition}

\newtheorem{cla}[prop]{Claim}

\newtheorem{remark}[prop]{Remark}

\usepackage{tikz}
\usetikzlibrary{arrows,decorations.pathmorphing,backgrounds,positioning,fit,petri}
\usetikzlibrary{shapes}
\usetikzlibrary{arrows.meta}

\newcommand{\tclock}[5]{
\begin{pgflowlevelscope}{\pgftransformscale{#4}}
\begin{scope}[shift={(#1,#2)}]
\shadedraw [inner color=#3!7!white, outer color=#3!90!black, 
    line width=0.2pt] (0,0) circle (0.5cm);
\foreach \x in {6,12,...,360} {\draw[line width=0.2pt] (\x:0.40cm) -- (\x:0.45cm);}
\foreach \y in {30,60,...,360} {\draw[line width=0.2pt] (\y:0.35cm) -- (\y:0.45cm);}
{\pgfsetarrowsstart{to}
\draw[line width=0.4pt] (0:0.29cm) -- (0.02,0);
\draw[line width=0.4pt]  (90:0.32cm)--(0,0.02);}
\filldraw[fill=black] (-0.055,0.55) rectangle (0.055,0.6);
\filldraw[fill=black] (-0.015,0.51) rectangle (0.015,0.55);
\draw [line width=0.2pt](0,0.61) circle (0.11cm);
\draw [line width=0.2pt](0,0) circle (0.5cm);
\draw [line width=0.2pt](0,0) circle (0.02cm);
\draw [red,thick,domain=30:45] plot ({#5*0.6*cos(\x)}, {#5*0.6*sin(\x)});
\draw [red,thick,domain=20:55] plot ({#5*0.65*cos(\x)}, {#5*0.65*sin(\x)});
\draw [red,thick,domain=10:65] plot ({#5*0.7*cos(\x)}, {#5*0.7*sin(\x)});
\draw [red,thick,domain=135:150] plot ({#5*0.6*cos(\x)}, {#5*0.6*sin(\x)});
\draw [red,thick,domain=125:160] plot ({#5*0.65*cos(\x)}, {#5*0.65*sin(\x)});
\draw [red,thick,domain=170:115] plot ({#5*0.7*cos(\x)}, {#5*0.7*sin(\x)});
\end{scope}
\end{pgflowlevelscope}
}

\newcommand{\Pb}{\mathbb{P}}
\newcommand{\E}{\mathbbm{E}}
\newcommand{\Id}{\mathbbm{1}}

\newcommand{\R}{\mathbb{R}}
\newcommand{\N}{\mathbb{N}}
\newcommand{\Z}{\mathbb{Z}}

\newcommand{\Ie}{\mathbbm{i}_\epsilon}

\title[Energy solutions to Stochastic Burgers equation]{On energy solutions to Stochastic Burgers equation}

\author{Patricia  Gon\c calves, \quad Alessandra  Occelli }
\address{Patricia Gon\c calves, Center for Mathematical Analysis,  Geometry and Dynamical Systems, Instituto Superior T\'ecnico, Universidade de Lisboa, Av. Rovisco Pais 1, 1049-001, Lisbon. E-mail: {\tt pgoncalves@tecnico.ulisboa.pt}}
\address{Alessandra Occelli, Center for Mathematical Analysis,  Geometry and Dynamical Systems, Instituto Superior T\'ecnico, Universidade de Lisboa, Av. Rovisco Pais 1, 1049-001, Lisbon. E-mail: {\tt alessandra.occelli@tecnico.ulisboa.pt}}
\thanks{P.G. thanks FCT/Portugal for support through the project UID/MAT/04459/2013. This project has received funding from the European Research Council (ERC) under  the European Union's Horizon 2020 research and innovative programme (grant agreement   No 715734). }
\makeindex
\begin{document}
\maketitle

\begin{abstract}
In this review we discuss  the weak KPZ universality conjecture for a  class of $1$-d systems whose dynamics conserves one or more quantities. As a prototype example for the former case, we will focus on weakly asymmetric simple exclusion processes, for which  the density is preserved and the equilibrium fluctuations are shown to cross from the Edwards-Wilkinson universality class to the KPZ universality class. The crossover depends on the strength of the asymmetry. For the latter case, we will present an exclusion process with three species of particles, known as the ABC model, for which we aim to prove the convergence to a system of coupled stochastic Burgers' equations, i.e. gradient versions of coupled KPZ equations. We will review the recent results on this matter around the notion of energy solutions introduced in~\cite{GJ1} for which uniqueness has been proved in~\cite{GP}.
	
\end{abstract}

\section{{Introduction}}
	\quad 
	
A classical problem in the field of interacting particle systems (IPSs)  is to derive the macroscopic laws of the
thermodynamical quantities of a physical system by considering an underlying microscopic dynamics which
is composed of particles that move according to some prescribed stochastic, or deterministic, law. The
macroscopic laws can be partial differential equations (PDEs) or stochastic PDEs depending on
whether one is looking at the convergence to the mean or to the fluctuations around that mean. The main goals of the field are related to the derivation of these macroscopic laws for general underlying dynamics. The IPSs were introduced in the mathematics community by Spitzer in~\cite{Spi} (but were already known to physicists) as microscopic stochastic systems, whose dynamics conserves a certain number of thermodynamical quantities of interest. The physical motivation for the study of IPSs is the following. Suppose that one is interested in analysing the evolution of some physical system, e.g. a fluid or a gas. Due
to a large number of molecules, it becomes hard to analyse the microscopic evolution of the system, so it is
more relevant to analyse the macroscopic evolution of its structure. Following Boltzmann's approach from
statistical physics, first one finds the equilibrium states of the system and then characterizes them through
macroscopic quantities, the thermodynamical quantities, such as: pressure, temperature, density, energy, etc. The
natural problem that arises then is to analyse the macroscopic behaviour of that system out of equilibrium.

The characterization and study of phenomena out of equilibrium is one of the biggest challenges of
statistical physics and despite its long history, a satisfactory answer to this kind of problems has not been
found yet~\cite{Spo}. From this approach, some PDEs arise that provide information about the macroscopic
evolution of the thermodynamical quantities of the system. Due to the huge complexity of the analysis of
these systems, some simplifications need to be introduced. For this purpose, one assumes that the underlying
microscopic dynamics, i.e. the dynamics between molecules, is stochastic, in such a way that a probabilistic
analysis of the system can be done. Assuming that the particles (or molecules) behave as interacting random
walks, subject to some random local restrictions, arise the so-called IPSs~\cite{Lig}.
Probabilistically speaking, these systems are continuous time Markov processes, being the time between
consecutive sites given by an exponential law. Due to the well-known property of memory loss of the exponential,
the future of the system, conditioned to its past, depends only on its present state.

In these notes we will be focused on  IPSs whose dynamics conserves some (one or more) thermodynamical quantities. 
By considering these microscopic random systems one can, by a space-time scaling limit, deduce (in a
rigorous way) the macroscopic laws governing the evolution of these quantities. These laws give the space-time
evolution of the thermodynamical quantities of the system and they are composed of one, or a system,
of PDEs. In the literature, this procedure is known as \emph{hydrodynamic limit}. The scaling limit is done by means
of a scaling parameter, that we denote by $n$, which allows us to identify relevant unknowns of the
macroscopic space with the respective microscopic space. Therefore if the macroscopic space is, for
example, the one-dimensional torus $[0,1)$, one can divide this space into $n$ subintervals of size $1/n$, in such a way
that each element $u\in[0,1)$  such that $u \in \Big[\frac{j}{n},\frac{j+1}{n}\Big)$, is identified with  $j$, namely with the closest integer number less or equal to $un$, so that the
microscopic space (a discrete space where the random dynamics is defined) is the set of points $\{0,1,..., n-1\}$,
the microscopic torus.

To summarize, on the underlying scenario one is interested in analysing, for example, the physical evolution of a gas  in a certain volume $V$. Two scales are considered: the macroscopic one, that describes the global motion of the gas; and the microscopic one, that describes the motion of the molecules, or particles, of the gas. One can think of these two scales in the following way. If we put a gas in a certain volume $V$ and with our eyes, we look at the evolution of the dispersion of the gas, this would be the macroscopic picture, while if we have a loupe to see how the molecules interact with each other, this would be the microscopic scale. Due to the large number of molecules, it is assumed that its motion is not deterministic but stochastic, in the sense that, each particle performs a random walk, that is, it moves, after some specified time, to a position of the discretized volume $V$. The discretization of the volume $V$ is done by means of the scaling parameter $n$. Under the assumption of the random movement of the molecules of the gas, a probabilistic analysis of the system can be performed. The process describing the random motion of this collection of particles is called an IPS. For details on the formal definition of an IPS we refer to~\cite{Lig}.

Initially,  on the volume $V$, we can  consider a point $u$ and a neighbourhood $V_u$ around that point. Since the number of particles is large, one can assume that locally the microscopic model is in equilibrium. The equilibrium is characterized by the thermodynamical quantity of interest, that we denote by $\rho(0,u)$. Now we let the system evolve in time. For  a time $t$, we pick the same point $u$ and assume again the local equilibrium, which is now characterized by $\rho(t,u)$. The important questions to answer are: \emph{How is the space-time evolution of $\rho(t,u)$? Is it described by some law? Which law is this?} 

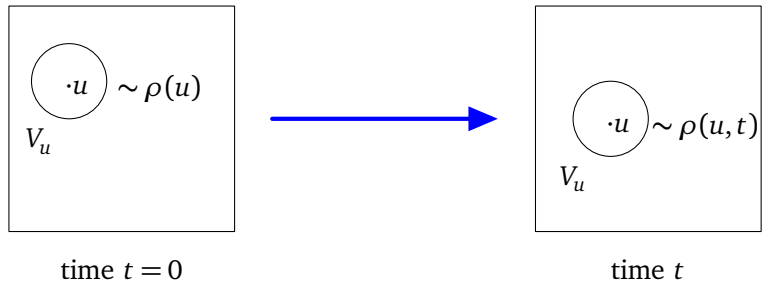
\begin{center}
	\tikzset{>=latex}
\begin{tikzpicture}[line cap=round,line join=round,>=triangle 45,x=1.0cm,y=1.0cm]
\draw [draw=black] (0,0) rectangle (3,3);
\draw [draw=black] (7,0) rectangle (10,3);
\node at (1.5,-.5) {{time $t=0$}};
\node at (8.5,-.5) {{time $t$}};

\draw [->,ultra thick, blue] (3.5,1.5) -- (6.5,1.5);
\draw[] (.8,2) circle (.5cm);
\draw[] (8,1.5) circle (.5cm);
\node at (.9,1.9) {$\cdot u$};
\node at (8.1,1.4) {$\cdot u$};
\node at (.4,1.2) {$V_u$};
\node at (7.5,.7) {$V_u$};
\node at (2,1.9) {$\sim \rho(u)$};
\node at (9.26,1.4) {$\sim\rho(u,t)$};
\end{tikzpicture}	
\end{center}

 So, the picture behind our interpretation is the following. To have a physical system evolving in a certain volume $V$, first one needs two scales: one for space and one for time. The evolution of its particles is random and Markovian. In the next step, one discretizes the volume $V$, according to the scaling parameter $n$, obtaining a certain number of cells. In each cell, one can put a random number of particles (for example, this number can be the result of some random experiment). At each cell, particles wait a \emph{random time} after which one of them decides to jump to some position of the discretized volume $V$ according to a given probability transition rate. The random times are assumed to be independent and exponentially distributed. In this case, the particle system falls into the class of \emph{Markov processes}, for which there is a vast and very well established theory.

	When analysing an IPS the first important task is to find its invariant measures. In many cases,  these measures are of product type and are known explicitly. It is not always an easy task to describe and characterize the invariant measures of general IPSs. Still nowadays, there are many IPSs whose invariant measures are not known explicitly, neither their structure nor their limiting properties. After this first task, one has to characterize these invariant measures by the thermodynamical quantities of the system such as the density, the pressure, the temperature, the energy, the viscosity, etc. The goal in the hydrodynamic limit theory consists in obtaining the macroscopic laws that govern the space-time evolution of each conserved quantity of the system. These macroscopic laws, in general PDEs, are called \emph{hydrodynamic equations}. 
 This result is a law of large numbers.
 
 The next step is to look at the deviations of the microscopic random system with respect to the deterministic limit. More precisely, we analyse the central limit theorem, namely the fluctuations of the conserved quantities of the system. In this respect, usually one looks at the deviations of the typical trajectory, which was obtained in the hydrodynamic limit, but starting the system from the invariant state or close to it. Nowadays it is still quite difficult to obtain the fluctuations starting from general measures since the estimate on space time correlation functions is not always easy and can be, in fact, quite challenging.

 In these notes we are going to analyse the fluctuations of some weakly asymmetric systems as exclusion processes, zero-range processes and Hamiltonian systems with more than one conserved quantity and we will show that the fluctuations will be governed by an Ornstein-Uhlenbeck process when the strength of the asymmetry is quite low, and by energy solutions of the Stochastic Burgers Equation, when the asymmetry is quite strong. 
 \\
 
 The outline of these notes is as follows. In Section~\ref{sec:hydro} we consider the one-dimensional exclusion process, a poster child to present the strategies and techniques to derive the hydrodynamic equations and the fluctuations. We will see that the corresponding (stochastic) PDEs depend on many factors, such as the scaling, the strength of the asymmetry, the length of the jumps. In Section~\ref{sec:weakKPZ}, in the framework of the weak KPZ universality conjecture, we show how to prove KPZ membership for a class of IPSs: their fluctuations are indeed described by the stochastic Burgers' equation, which is the gradient of the KPZ equation. In particular, we will see that these models exhibit a crossover behaviour from the Ornstein--Uhlenbeck process to the SBE, by tuning the asymmetry parameter.  In Sections~\ref{sec:tight} and~\ref{sec:BG} we prove tightness for the density field, by use of a second-order Boltzmann--Gibbs principle. A few examples of IPSs with one conservation law are provided in Section~\ref{sec:example}. Models with more than one conserved quantity are discussed in Section~\ref{sec:2CL}. In particular, in Section~\ref{sec:ABC} we focus on an IPS evolving  on the torus with two conserved quantities, known as  the ABC model.  We expect that, under a specific choice of the parameters tuning the asymmetry, the fluctuations are described by a system of coupled SBEs, or by independent SBE and OUE under further conditions.  This is justified by the predictions made using the nonlinear fluctuating hydrodynamics (NLFH) approach~\cite{Spo2} with a mode coupling treatment~\cite{SS}, explained in Section~\ref{sec:MCT}. This method allows to predict the universality class from the structure of the mode coupling
 matrices. These are characterized by the relation between the densities associated to the conserved quantities and the stationary currents. With this method we find the correct ``observables'' expressed as linear combinations of the  fluctuation fields of the conserved quantities.

\section{Hydrodynamics and fluctuations}\label{sec:hydro}

	\subsection{{The prototype  example: the exclusion process}} 
	We are going to state all the results  for the prototype example of an IPS, namely, the exclusion process  that we denote by $\eta(t)$, with $t$ running in  a compact set $[0,T]$, with $T$ fixed. To make the presentation as simple as possible we consider the process evolving on $\mathbb Z$, but other types of settings can be considered, such as $\Z^d$ or the torus $(\Z/N\Z)^d$. The dynamics of the simplest example of exclusion process is defined as follows. Each particle waits an exponential time of parameter $1$ and then it jumps to a site according to a certain probability transition rate. Due to the exclusion rule, the jump of a particle is performed if, and only if, the destination site is empty, otherwise nothing happens and the particle waits a new random time. The space state of this process is $\{0,1\}^{\mathbb Z}$. If another discrete set $\Lambda$ is considered, then the state space is $\{0,1\}^{\Lambda}$. Let us denote the jump rate from the site $x$ to the site $y$ by $p(x,y)$. In the figure below we represent a possible configuration with the description of the dynamics. 
	
 
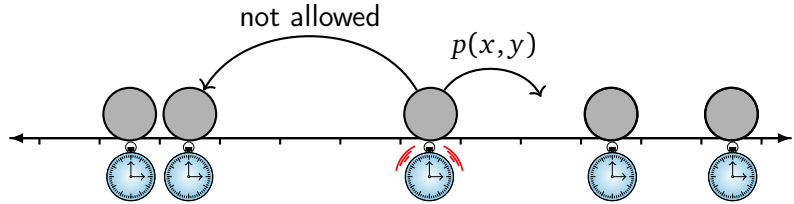
\begin{center}
 \begin{tikzpicture}[thick, scale=0.8][h!]
 \draw[latex-] (-6.5,0) -- (6.5,0) ;
\draw[-latex] (-6.5,0) -- (6.5,0) ;
\foreach \x in  {-6,-5,-4,-3,-2,-1,0,1,2,3,4,5,6}
\draw[shift={(\x,0)},color=black] (0pt,0pt) -- (0pt,-3pt) node[below] 
{};
 \node[fill=black!30!,shape=circle,draw=black,minimum size=0.7cm] (A) at (0.5,0.4) {};
    \node[fill=black!30!,shape=circle,draw=black,minimum size=0.7cm] (B) at (-4.5,0.4) {};
    \node[fill=black!30!,shape=circle,draw=black,minimum size=0.7cm] (C) at (-3.5,0.4) {};
        \node[fill=black!30!,shape=circle,draw=black,minimum size=0.7cm] (D) at (5.5,0.4) {};
   \node[fill=black!30!,shape=circle,draw=black,minimum size=0.7cm] (E) at (3.5,0.4) {};
    \node[shape=circle,draw=black,minimum size=0.7cm] (D) at (5.5,0.4) {};
    \node[shape=circle,draw=black,minimum size=0.7cm] (E) at (3.5,0.4) {};
        \node[shape=circle,minimum size=0.5cm] (G) at (2.5,0.4) {};
    \path [->] (A) edge[bend left=60] node[above] {$p(x,y)$} (G);
  \path [->] (A) edge[bend right=60,draw=] node[above] {\textsf{not allowed}} (C);
\tclock{0.6}{-0.8}{columbiablue}{0.8}{1}
\tclock{6.9}{-0.8}{columbiablue}{0.8}{0}
\tclock{4.4}{-0.8}{columbiablue}{0.8}{0}
\tclock{-4.4}{-0.8}{columbiablue}{0.8}{0}
\tclock{-5.6}{-0.8}{columbiablue}{0.8}{0}
 \end{tikzpicture}
 \end{center}

\quad

Here we assume that $p(x,y)=p(y-x)$, that is, the jump rate from the site $x$ to the site $y$ depends only on the distance between the departing and destination sites. 
	When jumps are allowed only to nearest neighbours  the process is said to be simple, so that $p(z)=0$ if $ |z|>1$. We consider 

\begin{equation}\label{eq:wek_as_rates}
p(1)=b_{+}+\frac {a}{n^\gamma}\quad  \textrm{and}\quad p(-1)=b_{-}-\frac {a}{ n^\gamma},\end{equation}

where $b_+,b_-, a$ and $\gamma$ are constants to be chosen ahead. 
If $a=0$ (there is no dependence on $\gamma$) and: 
\begin{enumerate}
	\item[I.] $b_{+}=b_{-}=1/2$, then the process is the symmetric simple exclusion process (SSEP);

\item[II.]  $b_{+} \neq b_{-}$ the process is the asymmetric simple exclusion process (ASEP).
\end{enumerate}

 If $a\neq 0$ (there is a dependence on $\gamma$) and $b_{+} = b_{-}$, the process is the weakly asymmetric simple exclusion process ($\gamma$-WASEP), being the jump rate given by a superposition of a symmetric dynamics by its asymmetric version. Note that the parameter $\gamma$ rules the strength of the asymmetry, the higher the value of $\gamma$ the weaker is the asymmetry. 
 
We note that we could also allow jumps not necessarily to nearest neighbour positions, like, for example, fixing a window of a certain size $K$ (as, for example,  in~\cite{DT}) and let particles evolve with similar rates, but this system has exactly the same asymptotic behaviour as the system with nearest neighbour jumps, see, for example,~\cite{BGS} and references therein. One could also  allow particles to jump to any site of the discrete lattice. In this case, we fix the probability transition rate to be, for example,  given by
\begin{equation}
p(z)=\frac{c_\pm}{|z|^{1+\alpha}},
\end{equation}
 with $\alpha>0$. The constants $c_\pm$ are defined in such a way that $p(\cdot)$ can be symmetric or not,  they can depend, for example, on the sign of the point $z$, see for example~\cite{GJ2}. The process obtained from this choice of $p(\cdot)$ is the exclusion process  with long jumps. 
 
 We note that the exclusion dynamics just defined does not destroy or create particles, since particles simply move in a discrete space according to some prescribed rule, therefore, the density of particles is a conserved quantity of the system and this is the thermodynamical quantity of interest when analysing this system.
 
 	\subsection{{Hydrodynamic limits}}
  Taking into account that the density of particles is the thermodynamical quantity of interest,   one defines a \emph{density empirical process}, which consists on a process of random measures, such that, for a fixed time $t$, each measure is a sum of Dirac measures supported at each site of $\mathbb{Z}$ that gives weight $\frac 1n$ to each particle, namely: 
\begin{equation}
\pi^n_t(\eta(t),du)=\frac{1}{n}\sum_{x\in\mathbb {Z}}\eta_x(t\theta(n))\delta_{x/n}(du).
\end{equation}
Above $\delta_{A}(du)$ is the Dirac measure of the set $A$.  Note that above since a re-scaling in space is done, one also needs to speed up the process in a time scale that we denote by $t\theta(n)$ and which will depend on the choice of the transition probability $p(\cdot)$.
	The goal in the hydrodynamic limit consists in showing that:
	
	\begin{enumerate}
	\item[I.] if we start the process $\{\eta(t\theta(n))\}_t$ from a collection of measures $\mu_n$ for which a law of large numbers holds, i.e.\,the sequence of random measures $\pi^n_0(\eta, du)$ converges, in probability with respect to $\mu_n$ and when $n$ is taken to infinity, to the deterministic measure $\rho(0,u)du$, where $\rho(0,u)$ is a real-valued measurable function; 
	
	\item[II.]  then, the same holds at later times $t$, that is, the random measure $\pi_t^n(\eta(t\theta(n)),du)$ converges, in probability with respect to $\mu_n(t)$, the distribution of $\eta(t\theta(n))$  and  when $n$ is taken to infinity, to the deterministic measure $\rho(t,u)du$, where $\rho(t,u)$ is the solution (usually in a weak sense) of a PDE, which is called, the \emph{hydrodynamic equation} of the system. 
	
	\end{enumerate}
	For the exclusion processes defined above, by a scaling limit procedure, one can get as hydrodynamic equations: 
	\begin{enumerate}
	\item[I.]  for the SSEP and re-scaling time diffusively $t\theta(n)=tn^2$, the heat equation given by
	\begin{equation}
\partial_t\rho(t,u)=\Delta \rho(t,u) \quad \textrm{[Heat eq.]}.
\end{equation}
\item[II.] for the $1$-WASEP (respectively for the ASEP) and by re-scaling time diffusively $t\theta(n)=tn^2$ (respectively, in the hyperbolic scale $t\theta(n)=tn$), the viscous Burgers equation (respectively, the inviscid Burgers equation),
\begin{equation}
\partial_t \rho(t,u)=\frac{1}{2}\Delta\rho(t,u)+ (1-2b_+)\nabla \Big(\rho (t,u)(1- \rho(t,u)\Big) \quad \textrm{[Viscous Burgers eq.]}
\end{equation}
\begin{equation}
\partial_t \rho(t,u)= (1-2b_ +)\nabla\Big(\rho(t,u)(1- \rho(t,u)\Big)  \quad \textrm{[Inviscid Burgers eq.]}.
\end{equation}

\item[III.] for the exclusion process with long jumps  and re-scaling time as $t\theta(n)=tn^\alpha$, the fractional Heat equation (when $c_{\pm}=c$) and the fractional Burgers equation (when $c_+>c_-$), 
\begin{equation}
\partial_t\rho(t,u)=-c(-\Delta^{\alpha/2})\rho(t,u), \quad \textrm{ [Fractional Heat eq.]}
\end{equation}
\begin{equation}
\partial_t \rho(t,u)= -\frac{c_++c_-}{2}(-\Delta^{\alpha /2}) \rho(t,u)+d(c_+,c_-)(-\nabla^{\alpha/2}) \rho(t,u)\quad \textrm{ [Fractional Burgers eq.]},
\end{equation}
where, above $\Delta$ and $\nabla$ denote, respectively, the Laplacian and the spatial derivative; $(-\Delta^{\alpha/2})$ and $-\nabla^{\alpha/2}$ denote, respectively,  the fractional Laplacian and the fractional derivative; $d(c_+,c_-)$ is a constant which depends on the underlying microscopic dynamics.
\end{enumerate}
For details we refer to, for example,~\cite{J1} and references therein.

One of the standard methods~\cite{KL} to prove the hydrodynamic limit consists in showing tightness of the sequence $\{\pi_t^n\}_{n}$ and then  characterizing the limit points of this sequence; this is the entropy method introduced by Guo, Papanicolau and Varadhan in~\cite{GPV}. Basically, tightness is a consequence of applying  the Aldous' criterion and  exploring some martingales in the context of Markov processes that can be associated to the density empirical measure. The most demanding and technical part is the characterization of limit points.  For this, one has to use the associated martingales and to prove some replacement lemmas which allow recognising the limit of the martingale as the weak solution to the corresponding hydrodynamic equation. The natural questions that arise after solving this issue are: 

\vspace{0.1cm}

$\bullet$ What are the fluctuations around the mean for each one of these models? 

$\bullet$ Is there a pattern of the macroscopic equations depending on the microscopic jump rates? 

$\bullet$ Are there universality classes which the models (with general features) belong to? 

$\bullet$ If so, what is the relationship between these universality classes?  

$\bullet$ Are these equations linked by some parameter depending on the underlying dynamics?

\subsection{{Fluctuations}} \label{sec:fluctuations}
In order to understand the questions raised in the previous subsection, first one considers an IPS starting from an equilibrium state, since the non-equilibrium scenario is much less understood and highly more complicated. In the case of the exclusion processes defined above, the invariant measures are the Bernoulli product measures, which are parametrized by a constant $\rho$ (that we denote here by $\nu_\rho$) and that are defined by:
	\begin{equation}
\nu_\rho\{\eta\in \{0,1\}^{\mathbb Z}: \eta_x=1\}=\rho.
\end{equation}
Then one defines, what is called the \emph{density fluctuation field}, which is  a linear functional defined on a function $f$, living in some proper space of test functions, as, for example, the Schwartz space $\mathcal S(\R)$, in the following way: first integrate $f$ with respect to the density  empirical  measure $\pi_t^n(\eta,du)$, then remove its mean (with respect to the invariant state) and then multiply it by $\sqrt n$ (so that we are in the Central Limit Theorem scaling), to get
\begin{equation}
\mathcal{Y}^n_t(f)=\frac{1}{\sqrt n}\sum_{x\in\mathbb{Z}}f(\tfrac{x}{n})[\eta_x({t\theta(n)})-\rho]. 
\end{equation}
 Above we subtracted $\rho$ from $\eta_x{(t\theta(n))}$ because  $\rho$ is the mean of $\eta_x{(t\theta(n))}$, with respect to the invariant measure $\nu_\rho$, and it does not depend on $x$: 
$
\mathbb{E}_{\nu_\rho}[\eta_x{(t\theta(n))}]=\rho.
$
 For IPSs where the invariant measure $\mu$ is not translation invariant, this mean will be defined as $\rho^n(x)=\mathbb{E}_{\mu}[\eta_x]$, where $\mu$ is the invariant measure. The same happens when one considers the system out of equilibrium, the only difference in that case is that the function $\rho_t^n(x)$ is given by $\mathbb{E}_\mu[\eta_x(t\theta(n))],$ where $\mu$ is the starting measure (not time invariant) and therefore $\rho_t^n(x)$ also depends on the parameter $t$. 

	We will focus these notes on an IPS starting from an equilibrium state (or very close to it) because even in this situation the theory is still far from being completely understood. 
	At this level, the goal consists in obtaining the stochastic PDEs ruling the evolution of the limiting process of the sequence $\{\mathcal{Y}^n_t\}_n$, denoted by $\mathcal {Y}_t$ (which is obtained from $\{\mathcal{Y}^n_t\}_n$ in a proper topology by sending $n$ to infinity). Depending on the prescribed dynamics one can derive, from different underlying stochastic dynamics, different stochastic PDEs using this procedure. This brings a physical motivation for the study of the stochastic PDEs that will emerge from the microscopic dynamics. For the exclusion processes introduced above, one can get:

\begin{enumerate}

\item[I.]  for the SSEP and rescaling time diffusively $t\theta(n)=tn^2$, the Ornstein--Uhlenbeck equation (OUE) given by
\begin{equation}
d\mathcal {Y}_t = A \Delta \mathcal{Y}_t dt + \sqrt C\nabla \mathcal W_t \quad \textrm{[OUE]}
\end{equation}
\item[II.] for the $\frac{1}{2}$-WASEP and by rescaling time diffusively $t\theta(n)=tn^2$ - the Kardar--Parisi--Zhang equation (KPZE)  (introduced in~\cite{KPZ}) or its companion, namely the stochastic Burgers equation (SBE), depending whether one is looking at the height fluctuation field or the density fluctuation field,
\begin{equation}
dh_t = A\Delta h_t dt + B(\nabla h_t)^2 dt + \sqrt C \mathcal {W}_t \quad \textrm{[KPZE]}
\end{equation}
\begin{equation}\label{eq:SBE}
d\mathcal{Y}_t = A \Delta\mathcal {Y}_t dt + B \nabla \mathcal {Y}_t^2 dt + \sqrt C\nabla \mathcal W_t \quad \textrm{[SBE]} 
\end{equation}
Note that for the choice $B=0$ the  SBE becomes an OUE, but for $B\neq 0$ there is a non-linear term in the SBE that requires a special treatment. 
\item[III.]  the exclusion process with long jumps and rescaling time as $t\theta(n)=tn^\alpha$, the fractional OUE (when $\alpha>3/2$) and the fractional SBE/KPZE  (when $\alpha=3/2$): 
\begin{equation}d \mathcal {Y}_t = A(\mathcal {L}_{\rho,\alpha} )^* \mathcal {Y}_t dt + \sqrt{C (-\mathcal L_{1/2})} \nabla \mathcal W_t \quad \textrm{[Fractional OUE]}
\end{equation}
\begin{equation}
d \mathcal{Y}_t =A (\mathcal{L}_{\rho,\alpha})^* \mathcal{Y}_t dt + B \nabla {\mathcal Y}_t^2 dt +  \sqrt{C (-\mathcal{L}_{1/2}) }\nabla \mathcal {W}_t \quad \textrm{[Fractional SBE]}.
\end{equation}
\end{enumerate}

Above, $A,B,C$ are constants, $\mathcal W_t$ is a white-noise, and $\mathcal{L}_{\rho,\alpha}$ is the generator of a $\alpha$-stable Lévy process, being
$(\mathcal{L}_{\rho,\alpha})^*$ its adjoint with respect to $\nu_\rho$ and $\mathcal L_{{1}/{2}}$ its symmetric part.

 In~\cite{GJ1, GJ3} it was proved that for a large class of weakly asymmetric simple exclusion processes (where $b_+=b_-=\frac 12$ and $a\neq 0$ above), depending on the range of the parameter $\gamma$, the density fluctuations cross from the OUE to the SBE.  More precisely, in a phase of weak asymmetry ($\gamma>1/2$), the density fluctuations are given by an OUE  (the same OUE  as in the SSEP, corresponding to the choice $a=0$ and $b_+=b_-=\frac 12$). This means that in this regime, the asymmetry is not seen at the macroscopic level. Nevertheless, in a phase of strong asymmetry ($\gamma=1/2$), the density fluctuations are given by energy solutions of the SBE. Therefore the system crosses from, what is called,  the Edward-Wilkinson universality class to the KPZ universality class, by changing the value of the parameter $\gamma$, which rules the strength of the asymmetry. 
	
	The latter results can be stated in terms of the KPZ equation by simply taking the current (height) fluctuation field instead of the density fluctuation field, see, for example,~\cite{GJ1}. The current process is defined in the following way. Let $x$ be a  site in $\mathbb Z$ and fix the bond $[x,x+1]$. The current process $J_x^n(t)$ at the bond $[x,x+1]$ counts the number of particles that jump from the site $x$ to the site $x+1$ minus the number of particles that jump from the site $x+1$ to the site $x$ during the time interval $[0,t\theta(n)]$. The current field is defined similarly to the density field but we replace $\eta_x{(t\theta(n))}$ by $J_x^n(t)$ and the corresponding means, see, for example,~\cite{GPS} for the exact definition. 
	
	The last result on the crossover was established for a large class of simple exclusion processes, with the restriction that the symmetric part of the dynamics is gradient and, roughly speaking, that there is a lower bound on the spectral gap of the dynamics restricted to a box of size $k$ given by $k^2$, see~\cite{GJ1, GJS} for the more precise statement on this bound. The gradient condition is defined as follows. Let $j^n_{x,x+1}(\eta)$ be the instantaneous current of the system at the bond $[x,x+1]$ for the configuration $\eta$. This current is defined as the difference between the jump rate from the site $x$ to the site $x+1$ and the jump rate from the site $x+1$ to the site $x$. We note that if $\mathcal L$ is the generator of the simple exclusion process, due to the conservation of the number of particles, $\mathcal L\eta_x({t\theta(n)})$ is written as the gradient of the instantaneous current, more precisely, 
\begin{equation}
\mathcal L\eta_x(t\theta(n))= j^n_{x-1,x} (\eta(t\theta(n)))- j^n_{x,x+1}(\eta({t\theta(n)})).
\end{equation}
 The gradient condition requires that the symmetric part of the current, which comes from the symmetric part of the generator $\mathcal L$, namely from 
$\mathcal S=\frac{\mathcal L+\mathcal L^*}{2},
$
 has to be written as the gradient of some local function. A local function is a function defined on the state space of the Markov process (in our setting $\{0,1\}^{\mathbb Z}$)  which depends on a configuration $\eta$ only through a finite number of its coordinates. We note that the gradient condition is restrictive, nevertheless, many systems do enjoy this property.

  	\section{{On the weak KPZ universality conjecture}}\label{sec:weakKPZ}
	
	In this section we discuss the aforementioned crossover from the OUE to the SBE/KPZE for some microscopic systems. 
	This is related to the  \emph{weak KPZ universality conjecture}, which states that a large class of $1$-d weakly asymmetric conservative systems should converge to the SB/KPZ equations. To prove the crossover for gradient systems, a possible strategy is to  follow the approach of~\cite{FGS, GJ1, GJ2, GJ3, GJS,GJSi,GPS} by showing tightness of the sequence of the density fluctuation fields $\{\mathcal Y_t^n\}_{n}$, and afterwards one has to characterize the limit  field. The latter  is the biggest difficulty that we one faces when proving this type of  crossover.

	To fix ideas let us consider the weakly asymmetric exclusion process $\eta(t)$ with rate given by \eqref{eq:wek_as_rates}, with a proper choice of {$b_+=1/2$} and {$a/2$} such that
	\begin{equation}\label{jump_rate_1}
p(1)=\frac 12+\frac {a}{2n^\gamma}\quad  \textrm{and}\quad p(-1)=\frac 12-\frac {a}{2 n^\gamma},\end{equation}
where {$a$} is a constant and $\gamma$ is a parameter ruling the strength of the asymmetry and belonging to the interval $[1/2,+\infty)$. We have to consider the system speeded up in the diffusive time scale, so that here {and in what follows } $\theta(n)=n^2$. This is the time scale in which we see non trivial scaling limits. 
The infinitesimal generator of the process  $\eta(t) $ is given on cylinder functions $f:\{0,1\}^\mathbb Z\to \mathbb R$ by
\begin{equation}\label{generator_wasep}
\mathcal L f(\eta)=n^2\sum_{x\in\mathbb R}c_{x,x+1}(\eta)(f(\eta^{x,x+1})-f(\eta)),
\end{equation} 
where $c_{x,x+1}(\eta)=p(1)\eta_x(1-\eta_{x+1})+p(-1)\eta_{x+1}(1-\eta_x)$ and $\eta^{x,x+1}$ is the configuration obtained from the configuration $\eta$ by swapping the occupation variables $\eta_x$ and $\eta_{x+1}$:
\begin{equation*}\eta^{x,x+1}(y)=\eta_{x+1}\textbf{1}_{y=x}+\eta_x\textbf{1}_{y=x+1}+\eta_y\textbf{1}_{y\neq x,x+1}.\end{equation*}
 We consider the system starting from the invariant state $\nu_\rho$, that is the Bernoulli product measure with parameter $\rho\in(0,1)$.   We associate to this system its density fluctuation field, defined  on  $f$ belonging to  $\mathcal S(\mathbb R)$, the Schwartz space,  as 
\begin{equation}
\mathcal{Y}^n_t(f)=\frac{1}{\sqrt n}\sum_{x\in\mathbb{Z}}f(\tfrac{x}{n})[\eta_x(tn^2)-\rho]. 
\end{equation}
In order to avoid complications we take $\rho=1/2$, otherwise, we would have to perform a Galilean transformation and re-centre the  field in a moving frame with velocity given in terms of $\rho-\tfrac 12$. 
Following the strategy explained in section \ref{sec:fluctuations}, in order to analyse the equilibrium fluctuations of this model, we need to show tightness of the sequence $\{\mathcal Y_t^n\}_n$ in the Skorohod space of càdlàg trajectories $D([0,T], \mathcal S'(\mathbb  R))$ with respect to the uniform topology and to characterize the limit point. The goal is to show that the limit $\mathcal Y$ solves some martingale problem which characterizes uniquely the solution. 
In~\cite{GJ1} it was proved that the limit field is a solution to an OUE in the case $\gamma>\tfrac 12$ and for $\gamma=\tfrac 12$ the fluctuations are given by the energy solution of the SBE.  Here we want to explain what is the difficulty that one faces when trying to prove this crossover from the OUE to the SBE. 

Since we are working with a continuous time Markov process, from Dynkin's formula, see for example  Lemma 5.1 of Appendix 1 of \cite{KL}, we can associate to the density fluctuation field a collection of martingales  given by
\begin{equation}\label{eq:dynkinWasep}
\mathcal M_t^n(f)=\mathcal Y_t^n(f)-\mathcal Y_0^n(f)-\int_{0}^t(\mathcal L+\partial_s)\mathcal Y_s^n(f)\, ds
\end{equation}
and with quadratic variation given by
\begin{equation}\label{quadvar}
\int_0^t \mathcal {L}  (\mathcal{Y}_s^n(f))^2 - 2 \mathcal  {Y}_s^n(f) \mathcal {L}\mathcal Y_s^n(f) ds.
\end{equation}
{These are martingales with respect to the natural  filtration $\mathcal F_t=\sigma(\mathcal Y_s(f), s\leq t, f\in \mathcal S(\mathbb R)).$} 
Up to very recently, there was not a proof of uniqueness of solutions of the SBE/KPZE. In the last few years, the first answer to this problem came from M. Hairer in~\cite{Hai} and~\cite{Hai2} with his revolutionary theory of regularizing structures. In those articles, he proposed a notion for solutions to the SBE/KPZE  and by smoothing the noise that appears in these equations, he could prove that those solutions are unique and that, in fact, they coincide with what is called \emph{ Cole-Hopf solutions}~\cite{BerG}. From the point of view of the IPSs, the imposed conditions on the solutions proposed by Hairer in~\cite{Hai} and~\cite{Hai2} are very difficult to check, so even after his incredible work, there was still a gap for showing the weak KPZ universality conjecture. In~\cite{GuJ} the authors proposed a notion of energy solutions to the (fractional) SBE which, in the case of the SBE,  is basically  the notion of energy solution of the SBE/KPZE proposed originally by~\cite{GJ1} and~\cite{GJS}. \\

\subsection{{Energy solution to the SBE}}\label{sec:SES}
In~\cite{GJ1}, the authors introduced the notion of energy solutions of the KPZE and SBE. They considered a one-dimensional, stationary, weakly asymmetric, conservative particle system {evolving on $\mathbb Z$}  and they proved that the limit point of the density fluctuations is given by energy solutions of the SBE. This notion of solution assumes the stationarity of the density field and that {$\mathcal Y_t\in\mathcal{S}'(\mathbb R)$} solves the martingale problem associated to \eqref{eq:SBE}, that is 
\begin{equation}
d\mathcal{Y}_t = A \Delta\mathcal {Y}_t dt + B \nabla \mathcal {A}_t^2 dt + \sqrt C\nabla \mathcal W_t,
\end{equation}
with {$\mathcal A_t=\lim_{\epsilon\to0}\mathcal{Y}_t\star i_\epsilon$}, where $i_\epsilon$ is an approximation of the identity.  This is a possible way to give sense to the nonlinearity in the equation above,  that allows to define a solution. However, to show some continuity properties of the quadratic term (which, in particular, implies that the quadratic variation is zero), another assumption is necessary. This is called ``energy condition'' or ``energy estimate''. {Let $\{\mathcal{Y}_t,t\in[0,T]\}$ be a stationary process. For $s<t\in[0,T]$, $\epsilon>0$ and $f\in\mathcal S(\R)$, let us define
	\begin{equation}
	\mathcal{A}^\epsilon_{s,t}(f)=\int_s^t\int_{\R}\mathcal{Y}_r(i_\epsilon(x))^2\nabla f dxdr.
	\end{equation} 
	We say that $\{\mathcal{Y}_t,t\in[0,T]\}$ satisfies an energy estimate if there exists a constant $\kappa$ such that
\begin{itemize}
	\item For any $f\in\mathcal S(\R)$ and any $s<t\in[0,T]$,
	\begin{equation}\label{eq:EE1}
	\E\Big[\Big(\int_s^t\mathcal{Y}_r(\Delta f)dr\Big)^2\Big]\leq\kappa(t-s)\|\nabla f\|_2^2;
	\end{equation}
	\item For any $f\in\mathcal S(\R)$, and $\delta<\epsilon\in(0,1)$ and any $s<t\in[0,T]$,
	\begin{equation} \label{eq:EE2}
	\E\left[\left(\mathcal{A}^\epsilon_{s,t}(f)-\mathcal{A}^\delta_{s,t}(f)\right)^2\right]\leq\kappa(t-s)\epsilon\|\nabla f\|_2^2.
	\end{equation}
\end{itemize}}

{In \cite{GJ1} it was introduced the following notion of energy solution.}
\begin{definition}~\label{def:energy_solution}
$\{\mathcal{Y}_t,t\in[0,T]\}$ is an energy solution to the SBE {given in \eqref{eq:SBE} } if
\begin{enumerate}
	\item[\bf{i)}] The process $\{\mathcal{Y}_t,t\in[0,T]\}$ is stationary and it satisfies the energy estimate;
	\item[\bf{ii)}] For any test function $f\in\mathcal S(\R)$ and any $t\in[0,T]$, the process
	\begin{equation}
	\mathcal{Y}_t(f)-\mathcal{Y}_0(f)-A\int_0^t\mathcal{Y}_s(\Delta f)ds -B\mathcal{A}_t(f)
	\end{equation}
	is a continuous martingale of quadratic variation $Ct\|\nabla f\|_2^2$.
	\end{enumerate}
	\end{definition}

We remark that the energy condition is strong enough to obtain some path properties of the
stationary energy solutions: the process $\{\mathcal{Y}_t,t\in[0,T]\}$ is a.s. \"Holder-continuous of index $\alpha$, for any $\alpha<1/4$.

The notions of energy solutions in~\cite{GJ1},~\cite{GJS} or~\cite{GuJ} use the classical formulation of {the martingale problem as above} and the difference between the notions is that~\cite{GuJ} requires, apparently, two extra conditions with respect to the notion given in~\cite{GJ1} and~\cite{GJS}, namely: 

\quad

{{\bf{a)}}} The process $\mathcal A_t$ corresponding to the non-linear term  $\nabla \mathcal Y_t^2/ (\nabla h_t)^2$ of the SBE/KPZE, has to be, almost surely, a process of null quadratic variation; 

\quad

{{\bf{b)}}} For a fixed time T, the reversed processes $\mathcal A_{T-t}$ and $\mathcal Y_{T-t}$ solve the same martingale equation as the original processes $\mathcal A_t$ and $\mathcal Y_t$, with respect to their own filtrations. 

\quad

All things considered, one can realize that the difference between the solutions is that~\cite{GuJ} requires only one extra condition with respect to~\cite{GJ1} and~\cite{GJS}, since the condition {{\textbf{a)}}}  above is a consequence of the notion introduced in~\cite{GJ1} and~\cite{GJS}, more precisely, of the energy condition, {see \eqref{eq:EE1} and \eqref{eq:EE2}}. We note that the energy estimate is a consequence of a second order Boltzmann-Gibbs principle which allows  replacing local functions of the Markov process $\{\eta({t\theta(n)})\}_{n}$ by functions of the conserved quantities of the system, namely the density of particles. There has been a lot of results on showing convergence to the energy solution of the SBE/KPZE from different  underlying microscopic dynamics and in~\cite{GP} the uniqueness of energy solutions was finally established. {The notion where the uniqueness holds is the one given in Definition \ref{def:energy_solution} and adding  the condition {\textbf{b)}} given above.} 

Among the aforementioned results, in~\cite{GJS} the crossover fluctuations for the density field are shown for a more general class of nearest-neighbour weakly asymmetric mass-conservative IPSs. They consider systems with a gradient type dynamics, whose invariant measures are not necessarily in a product form. WASEP, the zero-range process,
kinetically constrained exclusion systems and gradient exclusion systems with speed change, fall into this category. The result is formulated in the same language of energy solutions and  they were obtained thanks to a sharp estimate on the Boltzmann--Gibbs principle. But the main contribution of that work was to extend the proof for systems starting from initial distributions given by perturbations of the invariant measure under a certain ``bounded entropy'' condition.
	Very recently, \cite{yang2020} proposed an extension of the theory of energy solutions, showing density field crossover behaviour for particle systems starting out of stationary initial conditions and with unbounded number of particles per site. The author identifies the limit SPDEs for approximations of the density fluctuation field starting "almost" at equilibrium and then, due to the attractiveness of the system, shows that the approximations hold at any time.
 \subsection{Tightness}\label{sec:tight}
 With the goal of showing that the scaling limit of the $\gamma$-WASEP is an energy solution to the OUE/SBE, this section is dedicated to the proof of tightness for the process ${\mathcal Y}^{n}_t$  with respect to the Skorohod topology of $D([0,T];\mathcal{S}'(\R))$.
 By Mitoma's criterion \cite{mitoma}, it is enough to  prove tightness of the real-valued processes ${\mathcal Y}^{n}_t(f)$, for any $f\in\mathcal S(\R)$.
 To do that, we start decomposing the martingale~\eqref{eq:dynkinWasep}, with the generator given in \eqref{generator_wasep},  as 	
 \begin{equation}
 {\mathcal{M}}^n_t(f)={\mathcal Y}^{n}_t(f)-{\mathcal Y}^{n}_0(f)-{\mathcal I}^{n}_t(f)-{\mathcal K}^{n}_t(f),
 \end{equation}
 where
 \begin{equation}
 \begin{aligned}
 {\mathcal I}^{n}_t(f)&=\int_0^t ds {\mathcal Y}_s^n\left(~\tfrac 12\Delta_nf\right),\quad \quad
 {\mathcal K}^{n}_t(f)&=a \frac{\sqrt n}{n^\gamma}\int_0^t ds \sum_{x\in\mathbb{Z}}\nabla^+_nf\left(\tfrac{x}{n}\right)\bar{\eta}_x(sn^2)\bar{\eta}_{x+1}(sn^2).
 \end{aligned}
 \end{equation}
 Above $\Delta^+_n$ and $\nabla^+_n$ denote, respectively,  the discrete laplacian and the discrete derivative
\begin{equation}\label{eq:discretederivatives}
 	\begin{aligned}
 	\Delta_n^+ f&=n^2\left\{f\left(\tfrac{x+1}{n}\right)-2f\left(\tfrac{x}{n}\right)+f\left(\tfrac{x-1}{n}\right)\right\}\\
 	\nabla_n^+ f&=n\left\{f\left(\tfrac{x+1}{n}\right)-f\left(\tfrac{x}{n}\right)\right\}
 	\end{aligned}
 \end{equation}
 and $\bar \eta(x)$ denotes the centred variable $\eta(x)-\rho$. Recall that we are taking the starting measure the Bernoulli product measure $\nu_\rho$ with $\rho=1/2$.
 The quadratic variation  of the martingale in \eqref{quadvar} can be computed and it gives
\begin{equation*}
\int_0^t\frac{1}{2n}\sum_{x\in\mathbb{Z}}\Big(\nabla_n^+ f\Big(\tfrac{x}{n}\Big)\Big)^2\Big((\eta_x(sn^2)-\eta_{x+1}(sn^2))^2+\frac{{a}}{n^\gamma}(\eta_x(sn^2)-\eta_{x+1}(sn^2))\Big) ds.
\end{equation*}

 Therefore, $\lim_{n\to\infty}\mathbb{E}_{\nu}[(\mathcal{M}^n_t(f))^2]= t\chi(\rho)\|\nabla f\|_2^2$, where  $\chi(\rho)=\int_{\Omega}(\eta_x-\rho)^2\nu_\rho(d\eta)$ is the variance of $\eta_x$.

We can show tightness of the martingale sequence using Aldous' criterion, Chebyshev's inequality and straightforward computations on the quadratic variation of ${\mathcal M}^{n}_t(f)$. For any stopping time $\tau\leq T$,
\begin{equation*}
\begin{aligned}
&\Pb_\nu(|{\mathcal M}^{n}_{\tau+\sigma}(f)-{\mathcal M}^{n}_\tau(f)|>\epsilon)\leq\frac{1}{\epsilon^2}\E_\nu\left[({\mathcal M}^{n}_{\tau+\sigma}-{\mathcal M}^{n}_\tau)^2\right]\\&
\leq\frac{1}{\epsilon^2}\E_\nu\Bigg[\int_{\tau}^{\tau+\sigma}ds\frac{1}{2n}\sum_{x\in\mathbb{Z}}\Big(\nabla_n^+ f\Big(\tfrac{x}{n}\Big)\Big)^2\Big((\eta_x(sn^2)-\eta_{x+1}(sn^2))^2+\frac{{a}}{n^\gamma}(\eta_x(sn^2)-\eta_{x+1}(sn^2))\Big)\Bigg].
\end{aligned}
\end{equation*}

Therefore, last expectation is bounded from above by $\sigma( O(1)+O(n^{-\gamma}))$, which vanishes as $\sigma\to 0.$
  
 To prove tightness of the sequence $\{{\mathcal I}^{n}_t(f),t\in[0,T]\}_{n\in\N}$ we use the {Kolmogorov--Prohorov--Centsov's (KPC) criterion. Applying the Cauchy--Schwarz inequality twice, we see that
 \begin{equation*}
 \begin{aligned}
 \E_\nu\left[({\mathcal I}^{n}_{t}(f)-{\mathcal I}^{n}_s(f))^2\right]&=\E_\nu\Bigg[\Bigg(\int_s^t\mathcal Y_s^{n}(\Delta_n  f)\Bigg)^2\Bigg] \leq{{2\chi{(\rho)}}}{n}(t-s)^2\,\sum_{x\in\mathbb T_n}\left[\nabla_n^+f\left(\tfrac xn\right)\right]^2.
 \end{aligned}
 \end{equation*}

 This is enough to show that the sequence is tight with respect to
 the uniform topology of $\mathcal C([0,T];\R)$ and that any limit point has $\alpha$-H\"older-continuous trajectories for any $\alpha<1/2$.
 
 The last step is to prove tightness for the sequence $\{{\mathcal K}^{n}_{t}(f),t\in[0,T]\}$, using the same criterion.  This is a consequence of a second order Boltzmann--Gibbs principle, stated in Theorem~\ref{thm:BG} below. 
 As a consequence, one can show that there is a constant $C>0$ such that
 \begin{equation}\label{eq:tildeKtight}
 \begin{aligned}
 &\E_\nu\Bigg[\Bigg({\frac{~\sqrt n}{ n^\gamma}}\int_0^t\sum_{x\in{\mathbb {Z}}}\nabla^+_nf\left(\tfrac xn\right)\bar\eta_{x}(sn^2)\bar\eta_{x+1}(sn^2)ds\Bigg)^2\Bigg]\leq C{\frac{n}{n^{2\gamma}}}\left\{\frac{tL}{n}+\frac{t^2n}{L^2}\right\}\|\nabla^+_n f\|^2_{2,n}.
 \end{aligned}
 \end{equation}
 If we choose $L=\sqrt {t} n$,~\eqref{eq:tildeKtight} is bounded by $Ct^{3/2}\|\nabla^+_n f\|^2_{2,n}{n^{1-2\gamma}}$, satisfying the hypothesis of the KPC criterion.
 Above, for a function $f\in\ell^2(\Z)$ we have denoted
 \begin{equation}
 \|f\|^2_{2,n}:=\frac 1n\sum_{x\in\Z}f(x)^2<\infty.
 \end{equation}
 
 \subsection{The second order Boltzmann--Gibbs principle}\label{sec:BG}
 In this section we state and give a sketch of the proof of the second order order Boltzmann--Gibbs principle {presented in \cite{GJSi} and which does not require any knowledge on the spectral gap of the system, contrarily to the first proof of this principle presented in \cite{GJ1} that required the spectral gap of the system on a box of size $K$ to scale as $K^2$}. For $\ell\in\N$ and $x\in\Z$ we introduce the averages on boxes of size $\ell$, one to the right of $x$, one to the left:
 \begin{equation}\label{eq:boxes_averages}
 \ora{\eta}^\ell_x=\frac{1}{\ell}\sum_{y=x+1}^{x+\ell}\eta_y,\qquad \ola{\eta}^\ell_x=\frac{1}{\ell}\sum_{y=x-\ell}^{x-1}\eta_y.
 \end{equation}
 \begin{thm}\label{thm:BG}
 	There exists a constant $C=C(\rho)>0$ such that, for any $L\in\N$, any $t>0$, and for any function $v\in\ell^2(\Z)$
 	\begin{equation}\label{eq:BG}
 	\E_v\Bigg[\Bigg(\int_0^t\sum_{x\in\Z}v(x)\Big\{\bar\eta_{x}(sn^2)\bar\eta_{x+1}(sn^2)-(\ora{\eta}^{L}_x(sn^2))^2+\frac{\chi(\rho)}{L}\Big\}ds\Bigg)^2\Bigg]
 	\leq Ct\left[\frac{L}{n}+\frac{tn}{L^2}\right]\|v\|^2_{2,n}
 	\end{equation}
 \end{thm}
 A consequence of this principle is that, for $\gamma>1/2$, the field ${\mathcal{K}}^n_t$ vanishes as $n\to\infty$, while for $\gamma=1/2$, and for $L=\epsilon n$, it can be replaced by 
\begin{equation} \begin{aligned}
  &\int_0^t ds \sum_{x\in{\mathbb{Z}}}\nabla^+_nf\left(\tfrac{x}{n}\right)(\ora{\eta}^{\epsilon n}_x(sn^2))^2,
  \end{aligned}
  \end{equation}
  which is equal to 
  \begin{equation}
  \int_0^t ds \frac{1}{n}\sum_{x\in{\mathbb{Z}}}\nabla^+_n f\left(\tfrac{x}{n}\right)\Big({\mathcal Y}^n_s ( i_\epsilon(\tfrac{x}{n}))\Big)^2,
  \end{equation}
  leading to
  \begin{equation}
  \mathcal M_t^n(f)=\mathcal{Y}^n_t(f)-\mathcal{Y}^n_0(f)-\int_0^t ds {\mathcal Y}_s^n\Big(\frac 12 \Delta_nf\Big )-\int_0^t ds \frac{1}{n}\sum_{x\in{\mathbb{Z}}}\nabla^+_n f\left(\tfrac{x}{n}\right)\Big({\mathcal Y}^n_s (i_\epsilon(\tfrac{x}{n}))\Big)^2,
  \end{equation}
with {$i_\epsilon(x)$} an approximation of the identity, given by {$i_\epsilon(x)(u)=\frac{1}{\epsilon}\Id_{[x,x+\epsilon]}(u)$}. Observe that, {as we mentioned above,  for $\gamma>1/2$ the last display does not have the rightmost term since it vanishes when $n\to+\infty$, } and for that reason, in this regime of $\gamma$ the fluctuations are given by the OUE,  while for $\gamma=1/2$ that term survives and the fluctuations are given by the SBE.

 The idea of the proof of the previous theorem, consists in the following decomposition of the local function
 \begin{align}
 \bar\eta_{x}\bar\eta_{x+1}-(\ora{\eta}_x^{L})^2+\frac{\chi(\rho)}{L}
 =&\bar\eta_{x}(\bar\eta_{x+1}-\ora{\eta}^{\ell_0}_x)\label{eq:BG1}\\
 &+\ora{\eta}^{\ell_0}_x(\bar\eta_{x}-\ola{\eta}^{\ell_0}_x) \label{eq:BG2}\\
 &+\ola{\eta}^{\ell_0}_x(\ora{\eta}^{\ell_0}_x-\ora{\eta}^{L}_x)\label{eq:BG3}\\
 &+\ora{\eta}^{L}_x(\ola{\eta}^{\ell_0}_x-\bar{\eta}_x) \label{eq:BG4}\\
 &+\bar{\eta}_x\ora{\eta}^{L}_x-(\ora{\eta}^{L}_x)^2+\frac{(\bar\eta_{x}-\bar\eta_{x+1})^2}{2L} \label{eq:BG5}\\
 &-\frac{(\bar\eta_{x}-\bar\eta_{x+1})^2}{2L}+\frac{\chi(\rho)}{L}, \label{eq:BG6}
 \end{align}
 for $\ell_0\in\N$. The key ingredients lie in the estimates of the first four terms. To obtain a bound on {the term on the right hand side of} \eqref{eq:BG1},~\eqref{eq:BG2} and~\eqref{eq:BG4} we apply a ``one block estimate'' that allows to replace the local functions by functions of the density of particles in a box of finite size. We fix the size of the box $\ell_0\in\N$. Notice that $\ora{\eta}^{\ell_0}_x$ and $\ola{\eta}^{\ell_0}_x$ are local functions with zero mean with respect to the invariant measure and {whose} supports do not intersect. There exists a constant $C=C(\rho)$ such that, for any $t>0$ and for any $v\in\ell^2(\Z)$,\\
 \begin{equation}\label{eq:1BE1}
 \E_v\Bigg[\Bigg(\int_0^t\sum_{x\in\Z}v(x)\ola{\eta}^{\ell_0}(sn^2)\Big[\bar\eta_{x+1}(sn^2)-\ora{\eta}^{\ell_0}_x(sn^2)\Big]ds\Bigg)^2\Bigg]\leq C\frac{t\ell_0}{n}\|v\|^2_{2,n}.
 \end{equation}
 For~\eqref{eq:BG3} we need a multiscale argument to replace the average over a box of size $\ell_0$ with the average over a box of macroscopically small size $L$.	There exists a constant $C=C(\rho)$ such that, for any $\ell_0\leq L$, any $t>0$ and for any $v\in\ell^2(\Z)$,
 \begin{equation}\label{eq:MSA1}
 \E_v\Bigg[\Bigg(\int_0^t\sum_{x\in\Z}v(x)\ola\eta^{\ell_0}_x(sn^2))\Big[\ora\eta^{\ell_0}_{x}(sn^2)-\ora{\eta}^{L}_x(sn^2)\Big]ds\Bigg)^2\Bigg]\leq C\frac{t}{n}\left[L+\frac{\ell_0^2}{L}\right]\|v\|^2_{2,n}.
 \end{equation}
 The first step of the multiscale argument is the one-block estimate~\eqref{eq:1BE1}. In the second step we need to double the size of the box obtained in the one-block estimate. We fix $\ell_k\in\N$. Define $\ell_{k+1}=2\ell_k$. There exists a constant $C=C(\rho)$ such that, for any $t>0$ and for any $v\in\ell^2(\Z)$,
 \begin{equation}\label{eq:2B1}
 \E_v\Bigg[\Bigg(\int_0^t\sum_{x\in\Z}v(x)\ola{\eta}^{\ell_k}_x(sn^2))\Big[\ora\eta^{\ell_k}_{x}(sn^2)-\ora{\eta}^{\ell_{k+1}}_x(sn^2)\Big]ds\Bigg)^2\Bigg]\leq C\frac{t\ell_k}{n}\|v\|^2_{2,n}.
 \end{equation}
 We use this step a number of times which is logarithmic on the size of the box to obtain the estimate of the multiscale argument. Finally, we can replace a function of the
 density in a microscopic box of finite large size by a function of the density over a macroscopically small box.

 \section{Examples} \label{sec:example}
 In this section, we give a few examples of IPSs with their fluctuation results when one or more quantities are conserved.
 
  	\subsection{{Microscopic models with local defects}}
The previous crossover from the OUE to SBE obtained by tuning  the strength of the asymmetry,   can also be obtained for microscopic models with defects. By this we mean that one can take, for example, the exclusion process as defined above but the jump rate $p(\cdot)$ is perturbed at one particular bond, site or boundary, in such a way that it becomes more difficult to cross this bond, site or boundary, with respect to the others. The SSEP with a slow bond was studied in~\cite{FGN} and~\cite{FGN1} and its weakly asymmetric version was studied in~\cite{FGS}. The case of the exclusion process with a slow site was studied in~\cite{FGSu} and the exclusion process with the slow boundary was studied in~\cite{bmns}. The case with long jumps was studied in~\cite{BGJO, BGJO2,BGS,BGJS}. Let us explain what happens in the case of the SSEP with a slow bond. Suppose that the choice for the slow bond is $[-1,0]$, so that  and the jump rate at this bond is given by
\begin{equation}
p(1)= \frac{\alpha}{2n^\beta}+\frac{a}{2n^\gamma}\quad \textrm{ and}\quad  p(-1)=\frac{\alpha}{2n^\beta}-\frac{a}{2n^\gamma}
\end{equation}
and in all other bonds of $\mathbb Z$ it is given by \eqref{jump_rate_1}.
Above $\alpha, \beta, C,\gamma$ are real numbers which have some restrictions so that the dynamics is well defined. 
In order to have well-defined rates, we impose $\gamma\geq \beta$. Note that the higher the value of $\beta$, the weaker is the rate to cross the slow bond. In~\cite{FGS} it was proved that, for $\gamma<\beta<1/2$, the density fluctuations are given by the same OUE obtained in~\cite{FGN1}, where the symmetric version of the model is considered (that is $C=0$). But for $\gamma=1/2>\beta$ the fluctuations are given by the SBE/KPZE.

 In~\cite{FGN},~\cite{FGN1} and~\cite{FGS} the symmetric and weakly asymmetric versions of this model were considered and the results showed a phase transition at the level of density of particles, current and tagged particle. More precisely: 

\quad

{\bf{1)}} for a phase of strong strength of the slow bond the asymptotic behaviour is the same as if there was no slow bond, so one gets: the heat equation with no boundary conditions as hydrodynamic equation, an OUE for the density fluctuations and fractional Brownian motions for the limiting distributions of the current and tagged particle; 

\quad

{\bf{2)}} for a phase of weak strength of the slow bond, one gets: the heat equation with boundary conditions of Robin type, stating that the amount of flux through the slow bond is proportional to the difference of concentration in the two regions separated by the slow bond (the Fick's law); the fluctuations are given by an OUE with Robin boundary conditions and the limiting distributions of current and tagged particle are given by a Gaussian process which is not self-similar, in particular, it is not a fractional Brownian motion;

\quad

{\bf{3)}} for a phase of a very weak strength of the slow bond, one gets: the heat equation with Neumann's boundary conditions as hydrodynamic equation, the OUE ruling the density fluctuations, but the limiting distributions of current and tagged particle are again given by fractional Brownian motions as it happens in  the case {\bf{1)}}. 
 	\subsection{{Microscopic models in contact with reservoirs}}
	One can also look at microscopic models in contact with reservoirs, in order to obtain the crossover from the OUE to the SBE/KPZE, but with boundary conditions. In case of the exclusion process, the definition of the dynamics can be seen for example in~\cite{GLM}, but the idea is simple. 
	For example, the SSEP in contact with one type of stochastic reservoirs as in~\cite{GLM} is defined in the following way. Consider,  the set of points $\{1,...,n-1\}$ where particles evolve under the exclusion rule and  add the sites $x=0$ and $x=n$ which will act as reservoirs. Suppose  that those sites $x=0$ and $x=n$ inject and remove particles at the sites $x=1$ and $x=n-1$ respectively.  Particles placed in $\{1,...,n-1\}$ jump to one of the nearest-neighbours at rate $1/2$ and the injection and removal rate of particle at the boundary is given by: the injection at $x=1$ is $\alpha/n^\theta$ and at $x=n-1$ is  $\beta/n^\theta$ and the removal rate at $x=1$ is $(1-\alpha)/n^\theta$ and at $x=n-1$ is $(1-\beta)/n^\theta$. The parameters $\alpha, \beta\in(0,1)$ and $\theta\in\mathbb R$. Depending on the value of $\theta$ the asymptotic behaviour of the system changes. For some examples on other types of exclusion processes in contact with reservoirs we refer the interested reader to the survey lecture notes~\cite{G1} and references therein. 
	
	We note that for this dynamics, the number of particles is no longer conserved because particles can enter or leave the system through the reservoirs. By imposing different densities at the reservoirs (for example, $\alpha$ at the left reservoir and $\beta$ at the right reservoir, with  $\alpha\neq\beta$) there is a natural flux of particles due to the difference of concentration of particles in the reservoirs,  but there is also another flux, which can come from the dynamics in the bulk, if the rates are not symmetric. For example, if the left reservoir has more particles than the one at the right, then particles will more often enter to the system from the left boundary and exit the system from the right boundary. If the rate in the bulk gives a preference to move to the right, then both fluxes cooperate. Nevertheless, if the rate in the bulk gives preference to move to the left, then there is a competition between the fluxes. 
	
When looking at the fluctuations of exclusion processes in contact with reservoirs one interesting problem is to analyse the weakly asymmetric model and to prove the crossover from the OUE to the SBE/KPZE with different boundary conditions. In~\cite{GPS} the case with Dirichlet boundary conditions was derived but for other types of equations and other boundary conditions the problem is open and it is quite interesting. 

 One alternative approach to this last problem is to consider the microscopic version of the Cole-Hopf transformation due to Gartner~\cite{Gar}, and used in~\cite{DG,BerG} and later in~\cite{DT,GLM,Corwin2016}, which is obtained by taking the exponential of the current process associated to the underlying microscopic dynamics. The advantage of this new process is that it solves a linear equation, similarly to what happens at the macroscopic level with the Cole-Hopf solution of the KPZE. Nevertheless, one has to be very careful when redefining the microscopic Cole-Hopf transformation, since it solves a linear equation (the non-linearity is removed by the exponentiation of the process) but with non-trivial non-linear boundary conditions. For a particular case of the dynamics described above, we refer the reader to~\cite{Corwin2016} and references therein. 
   
\section{{Models with several conservation laws}}\label{sec:2CL}

	One can also analyse microscopic models with more than one conserved quantity. In~\cite{BFS} multi-component coupled SBEs have been obtained as scaling limits of the empirical measures of a multi-species zero-range process. Here the method is still based on martingale representations of the fluctuation fields. But this time to write Dynkin's formula in a closed form, one needs a multi-component Boltzmann-Gibbs principle, which allows to replace nonlinear terms with functionals of the fields. Moreover, also the notion of energy solution is extended from a scalar to a multi-component formulation. Uniqueness is a consequence of~\cite{GP,GP2}.
			
	In~\cite{BS} it was introduced and studied numerically, a class of Hamiltonian models perturbed by conservative noises which present strong analogies with the standard chains of oscillators.
 Since the aforementioned  models have more than one conserved quantity (the particle densities for the multi-species zero-range, the energy and the volume for the Hamiltonian model), the analysis of its asymptotic behavior is more intricate than for the case of models with just one conserved quantity, like exclusion processes defined above. As examples of the Hamiltonian systems introduced above are the models of~\cite{BG,BGJ,BGJSS}. In those articles, the deterministic dynamics was perturbed by a weak strength noise and it was studied the fluctuations of the conserved quantities, namely, the energy and the volume, starting the system from the invariant measure, which is of product form. One can consider the weakly asymmetric version of those models (by increasing the strength of the noise) and analyse their crossover fluctuations. 
	
	The first problem in these models is the question of which fields to look at. Another problem  when studying these systems is that the conserved quantities can be linearly transported in the system, each one having its own velocity and its own time scale. When writing down the equations for the fields associated to the conserved quantities of the system, one has to remove the velocity in each field. Afterwards, one has to close the equations of each one of these fields, but each one evolves in a certain time scale (which is not necessarily the same for both) and one has to analyse the leading terms in the expansion of the equations in such a way that one can recover the limiting  stochastic PDE. 

	The main difficulty that we will face when looking at the field of the conserved quantities of these models is when one has to recognise the stochastic PDE satisfied by each limit field. At that point, one has to write the field associated to the instantaneous current of the system in terms of the field of the conserved quantities of the system,  that is, one has to derive a Boltzmann-Gibbs principle, as we did in Section~\ref{sec:BG} for one conserved field.

	The Boltzmann-Gibbs  principle we introduced was formulated for systems with one conserved quantity and it states that any local field of the dynamics can be replaced (in a proper topology) by the fluctuation field of the conserved quantity. When one is looking at the field of the conserved quantity of IPSs with asymmetric jump rates in the hyperbolic time scale $tn$ or symmetric jump rates in the diffusive time scale $tn^2$, the aforementioned Boltzmann-Gibbs principle is sufficient to recover the  stochastic PDE satisfied by the limiting field.  But when one is looking at an IPS in a reference frame moving at a certain velocity and in a time scale which is longer than the hydrodynamical one, then the Boltzmann-Gibbs principle of~\cite{BR} does not give any information about the limiting field. For that purpose, a second-order Boltzmann-Gibbs principle has been derived in~\cite{G,GJ1,GJS}. For exclusion processes, since the density is preserved, that principle allows  replacing certain local fields of the dynamics by the square of the density fluctuation field. After using this principle it becomes simple to close the equation for the field of the conserved quantity and to recognise the stochastic PDE satisfied by the limit field. 
		
	In~\cite{GJSi} the proof for the second order Boltzmann-Gibbs principle  does not impose strong conditions on the underlying microscopic dynamics and allows one to obtain the crossover fluctuations from the OUE to SBE/KPZE as described above for the $\gamma$-WASEP and, more generally, for any system which has a current written as a sum of polynomial functions. We believe that the extension of this result to Hamiltonian systems and to multi-species models~\cite{BFS,GO} will open a new way to obtain the fluctuations of the conserved fields and to establish the precise dependence on the strength of the perturbing noise, at the level of the crossover from different  stochastic PDEs. 
	
In the next section we consider a model with more than one conserved quantity and we  exemplify (as we did earlier for the $\gamma$-WASEP) how to derive the fluctuations of certain quantities of the system and how to overcome the problems that we  just mentioned above, namely, the choice of the fields and closing the equations.
	
\subsection{A prototype model: the ABC model with external field}\label{sec:ABC}	
Here we consider a system of particles of three species $\alpha\in\{A,B,C\}$, evolving by nearest neighbour particles exchanges in the discrete ring with $n$ sites, $\mathbb{T}_n=\mathbb Z/n\mathbb Z$, in the presence of an external field that interacts with particles of different species with a different strength; or, equivalently, the system ``sees'' three fixed constant fields $E_A,E_B,E_C$.
The total number of particles of each species $n_{\alpha}$, $\alpha\in\{A,\,B,\,C\}$, is conserved and it holds $n_{A}+n_{B}+n_{C}=n$.
The space of configurations is $\tilde{\Omega}_{n}=\{A,\,B,\,C\}^{\mathbb{T}_n}$; its elements are denoted by $\eta$; on each site $x\in\mathbb{T}_n$, $\eta_x\in\{A,\,B,\,C\}$.
The ABC model is a continuous time Markov chain with space state ${\Omega}_{n}$ and generator $\mathcal L_{n}$ acting on  functions 
$f\,:{\Omega}_{n}\rightarrow\mathbb{R}$ as

\[
\mathcal{L}_{n}f(\eta)=n^{{2}}\sum_{x\in\mathbb{T}_n}c^E_{x}(\eta)[f(\eta^{x,x+1})-f(\eta)].
\]

The rates $c^{\alpha\beta}_x(\eta)$ are given in the \emph{weakly asymmetric regime}: a configuration $(\alpha,\beta)$ on the bond $x,x+1$ is exchanged to $(\zeta,\xi)$ with rate $c^{\alpha\beta}_{x}=e{\frac{E_\alpha-E_\beta}{2n^\gamma}}$, for $\alpha,\beta\in\{A,B,C\}$. This is a generalization of the model originally introduced in~\cite{evans,evans2}. The stationary distribution $\nu$ is a product measure over $x\in\Z_N$ with parameters $\rho_A$, $\rho_B$ and $1-\rho_A-\rho_B$ corresponding to each of the states a site can assume, species $A$, $B$ or $C$/empty.

\subsection{{Hydrodynamics of the ABC model}}
We define the occupation numbers of the species $\alpha$  as the functions $\xi^{\alpha}:{\Omega}_{n}\rightarrow\{0,\,1\}^{\mathbb{T}_n}$ acting on the configurations in the following way:

\begin{equation}
\xi^{\alpha}_x(\eta):={\mathbbm{1}}_{\{\alpha\}}(\eta_x)=\begin{cases}
1 & \eta_x=\alpha\\
0 & \textrm{otherwise}.
\end{cases}
\end{equation}

Since it holds $\sum_{\alpha}\xi^{\alpha}_x(\eta)=1$
$\forall x\in\mathbb{T}_n$ and $\sum_{x\in\mathbb{T}_n}\xi^{\alpha}_x(\eta)=n_{\alpha}$, we can equivalently see the model as a two-species particle system, $A$ and $B$, and analyse the law of large numbers for $\xi^A_x$ and $\xi^B_x$.
In the diffusive time scaling $a=2$ and under $\gamma=1$, the hydrodynamic equations (\cite{Bertini} using the methods of~\cite{KL,KOV}) for the densities of particles $A$ and $B$ are given by 
\begin{equation}\label{eq:hydroeq_E}
\begin{aligned}
\partial_t \rho^A&=\partial_{xx}^2 \rho^A-\partial_x [\rho^A(1-\rho^A)(E_A-E_C)-\rho^A\rho^B(E_B-E_C)]\\
\partial_t \rho^B&=\partial_{xx}^2 \rho^B-\partial_x [\rho^B(1-\rho^B)(E_B-E_C)-\rho^A\rho^B(E_A-E_C)],
\end{aligned}
\end{equation}
and  the density of particles $C$ can be recovered as $\rho^C=1-\rho^A-\rho^B$. Observe that~\eqref{eq:hydroeq_E} can be written as
\begin{equation}
\partial_t \begin{pmatrix}\rho^A\\ \rho^B \end{pmatrix}=\Delta \begin{pmatrix}\rho^A\\ \rho^B \end{pmatrix}-\nabla\cdot \left(\rm{X}(\rho)\cdot g_E \right),
\end{equation}
where \begin{equation}
\rm{X}(\rho)=\begin{pmatrix} \rho^A(1-\rho^A) & -\rho^A\rho^B\\ -\rho^A\rho^B & \rho^B(1-\rho^B) \end{pmatrix}
\end{equation} is the mobility and $g_E=\begin{pmatrix}E_A-E_C\\E_B-E_C\end{pmatrix}$.

\subsection{{Fluctuations of the ABC model}}
Now we want to characterize a central limit theorem. The fluctuation results for this model have been obtained by the authors in an upcoming paper~\cite{GO}. {Since our dynamics now is evolving on $\mathbb T_n$ we need to take a suitable space for the fluctuation fields. We denote the space of smooth $\mathbb R$-valued
functions $f$ defined on $\mathbb T$ by ${\mathcal D(\mathbb T)}$ and we denote by ${\mathcal D'(\mathbb T)}$ the space of $\mathbb R$-valued distributions on $\mathbb T$. We define the density fluctuation fields ${\mathcal Y}_t^N=({\mathcal Y}_t^{N,A},{\mathcal Y}_t^{N,B},{\mathcal Y}_t^{N,C})\in {\mathcal D'(\mathbb T)}$ as } 
\begin{equation}\label{eq:densityfieldNL}
{\mathcal Y}_t^{N,\alpha}(du)=\frac{1}{\sqrt n}\sum_{x\in\mathbb{T}_n}\left[\xi^\alpha_x(tn^2)-{\mathbbm{E}}_\nu[\xi^\alpha_x(tn^2)\right]T^+_{v_\alpha(n) t}\delta_{\frac{x}{n}}(du),
\end{equation}
where ${\mathbbm{E}}_\nu[\cdot]$ is the expectation with respect to the invariant measure and $T^+_{v_\alpha(n) t}$ is the translation operator acting on $f\in{\mathcal D(\mathbb T)}$ as
\begin{equation}
T^+_{v_\alpha(n) t}f\left(\tfrac{x}{n}\right)=f\left(\tfrac{x+v_\alpha(n) t}{n}\right).
\end{equation}
{As we will see later in  Section~\ref{sec:convSBE}, the choice for the velocity  $v_\alpha(n)$ will depend on the parameter $\gamma$ and the role of this quantity is the centring of the density field in a proper moving frame.  Observe that $\mathcal L_n \xi_x^\alpha=j^\alpha_{x-1,x}-j^\alpha_{x,x+1}$, where $j_{x,+1}^\alpha$ is the the infinitesimal current for the species $\alpha$.  A simple computation shows that 
\begin{equation}
\begin{aligned}
j^\alpha_{x,x+1}=&\xi^{\alpha}_{x+1}-\xi^{\alpha}_{x}\\
&+\frac{E_{\alpha+1}-E_{\alpha}}{2n^\gamma}(\xi^{\alpha}_{x}\xi^{\alpha+1}_{x+1}+\xi^{\alpha+1}_{x}\xi^{\alpha}_{x+1})\\
&+\frac{E_{\alpha+2}-E_\alpha}{2n^\gamma}(\xi^{\alpha}_{x+1}+\xi^{\alpha}_{x}-2\xi^{\alpha}_{x}\xi^{\alpha}_{x+1}-\xi^{\alpha}_{x}\xi^{\alpha+1}_{x+1}-\xi^{\alpha+1}_{x}\xi^{\alpha}_{x+1}).
\end{aligned}
\end{equation}
The infinitesimal current for particles of species $\alpha+1$ is given by the same expression with $\alpha$ and $\alpha+1$ interchanged.
{Since $\mathcal L_n \xi_x^\alpha$ is given by a gradient, and the infinitesimal current  has a gradient term (the first term on the right hand side of last display)  which does not depend on $\gamma$, this term will give the contribution to the Laplacian in the equation. }
Observe that only the last two terms on the right hand-side of the last display depend on $\gamma$ and the current for species $\alpha$ involves terms of the species $\alpha+1$. Moreover, if   we scale the asymmetry as $1/n^\gamma$ with $\gamma>1/2$, then the only surviving terms are the ones giving the Laplacian  and the equations of the density field for the species $\alpha$ is independent of the species $\alpha+1$.  The same  happens when we look at the current for the species $\alpha+1$. Thus, we would have a diffusive behaviour for both density fields. Nevertheless, for  $\gamma=1/2$, the second and third terms on the right hand-side of the last display can contribute to the limit and in that case we can  try to observe KPZ behaviour. In the next section we explain how to deduce what are the fields that one should look at, when dealing with systems with more than one conservation law and we apply the results to the ABC model. 

\subsection{Mode coupling theory}\label{sec:MCT}
	The description of the universality classes from OUE to SB/KPZ equation is for one component systems, that is, for systems conserving only one quantity. In the case of multicomponent systems, as the Hamiltonian systems or the ABC model just mentioned, the scenario is much less understood. Nevertheless, this problem has recently attracted a lot of interest in both the physics and the mathematics community. This problem appeared, in particular, in the context of superdiffusion of energy in $1$-d nonlinear chains of oscillators. The high interest of physicists on the superdiffusion problem can be understood by the fact that, up to very recently, no theoretical explanation was available and more puzzling is that the numerical simulations are too controversial to provide any clear insight, see, for example,~\cite{LLP,D}.
	
	Recently,~\cite{Spo2,SS} with the so-called "nonlinear fluctuating hydrodynamics theory", which has been developed during the last three years, proposed a rich and complex phase diagram of the universality classes for these Hamiltonian systems. The richness of the diagram is explained by the non-trivial nonlinear couplings, occurring at different time scales, between the conserved quantities. This means, as mentioned above, that each conserved quantity will have its own velocity in a certain time scale and its own limiting  stochastic PDE ruling its fluctuations. The only non-KPZ universality class which has been, so far, proved rigorously to occur is in the context of harmonic chains perturbed by a conservative noise (see~\cite{JKO,BGJ}) or with a very small anharmonicity (see~\cite{BGJS}). 
	
	The predictions in~\cite{Spo2,SS} are done in the context of chains of oscillators by using the following idea. First, they linearize the hydrodynamic equation of a general system, which is supposed to be a Euler equation given by 
	\begin{equation*}
	\partial_t u'+\partial_x\langle g'\rangle_{u'}=0,
	\end{equation*}
	where $u'=(u_1,u_2,u_3)$ is the vector whose $i$-th $(i=1,2,3)$ component $u_i$ represents one of the conserved quantities of the system, $g'$ is the vector whose $i$-th component $g_i$  represents the instantaneous flux of the system for one of the conserved quantities and $\langle\cdot\rangle_{u'}$ represents the average with respect to the invariant state with parameter given by $u'$, that is, the local Gibbs state. By linearizing the equation above, they get to $\partial_tu'+A\partial_xu'=0$. Adding a noise term and a dissipating term, one arrives at the following equation 
	\begin{equation*}
	\partial_tu'+\partial_x[Au'-D\partial_xu'+BW']=0
	\end{equation*}
	where $W'$ is an $n$-dimensional white noise and $A$, $D$ are matrices. Considering the dynamical correlation function $S(j,t)$ given by 
$
	S(j,t)=\langle g_j(t)g_0(0)\rangle
	$
	where the average is taken with respect to the invariant measure, then $\sum_j jS(j,t)=ACt$, where $C$ is another matrix satisfying $DC+CD=BB^T$. The predictions are made by taking certain ansatz for the matrices given above and they are quite impressive, but, on the other hand, they are not mathematically rigorous. Nevertheless, they bring up a new insight to approach the problem from the mathematical point of view. It is clear now that in order to study the fluctuations of systems with multicomponent conserved quantities one has to look at a proper linear combination of the fields of the conserved quantities. Up to now, it was not even possible to pose this problem in a formal mathematical way, but according to these predictions, this is exactly what one has to do. 
	
	The articles~\cite{Spo2} and~\cite{SS} gave very detailed predictions about the correct time scales that one should see a non-trivial behavior for each one of the conserved quantities of the system, and more than that, they also predicted what are the limiting processes that one is searching for. They did it for the models of chains of oscillators with more than two conserved quantities, but these models should have the same asymptotic behavior as the dynamics introduced above with only two conserved quantities. According to their predictions, one can get crossovers from OUE to KPZE, KPZ fixed point, and Lévy processes with an exponent which can be any number of the Fibonacci sequence, for example. There is a limiting Lévy process that is called the gold Lévy process, that is a Lévy process of exponent equal to the gold number $\frac{1+\sqrt 5}{2}$. 

\subsection{{Predictions from MCT}}\label{sec:modeABC}

Since the system has two conserved quantities, the densities $\rho^A$ and $\rho^B$, it is possible to have an estimate on the scaling exponent with mode coupling theory, following the approach of~\cite{Popkov2015,Spo2,SS}. Without loss of generality one can consider the case $E_{\alpha+2}=0$. If one wants to recover the original model, one has just to replace $E_{\alpha(+1)}$ with $E_{\alpha(+1)}-E_{\alpha+2}$.
To obtain the normal modes of the system  one needs to find the matrix $R$ that diagonalizes the Jacobian $J$ of the current, $j(\rho)=-\nabla\rho+\chi(\rho)E$. The eigenvectors and the eigenvalues of $J$ determine the observables to study and velocity of the moving frame for each field.
In the equal density case $\rho_\alpha=\rho_{\alpha+1}=1/3$, the normal modes are given by
\begin{equation}
\begin{aligned}
\mathcal Z_t^n&=E_\alpha\mathcal Y^{n,\alpha}_t+[E_\alpha-E_{\alpha+1}-\delta]\mathcal Y^{n,\alpha+1}_t,\\
\widetilde{\mathcal Z}_t^n&=E_\alpha\mathcal Y^{n,\alpha}_t+[E_\alpha-E_{\alpha+1}+\delta]\mathcal Y^{n,\alpha+1}_t,
\end{aligned}
\end{equation}
with
\begin{equation}\delta=\frac 23\sqrt{E_{\alpha}^2+E_{\alpha+1}^2-E_{\alpha}E_{\alpha+1}}.\end{equation}
According to the predictions, under the choice $a=2$ and $\gamma=1/2$, both modes are supposed to be KPZ. 

Though, if we assume certain hypothesis  on the rates, there are two cases for which MCT give different predictions:
\begin{enumerate}
	\item[I.] when $E_{\alpha+1}-E_{\alpha}=E_{\alpha+2}-E_{\alpha}=E$, the normal modes are 
	\begin{equation}\label{eq:modes2}
	\mathcal Z_t^n=\mathcal Y^{n,\alpha}_t,\quad \widetilde{\mathcal Z}_t^n=\mathcal Y^{n,\alpha}_t+2\mathcal Y^{n,\alpha+1}_t;
	\end{equation}
	\item[II.] when $E_{\alpha}-E_{\alpha+2}=E_{\alpha+1}-E_{\alpha+2}=E$, the normal modes are
	\begin{equation}
	\mathcal Z_t^n=\mathcal Y^{n,\alpha}_t+\mathcal Y^{n,\alpha+1}_t,\quad \widetilde{\mathcal Z}_t^n=\mathcal Y^{n,\alpha}_t-\mathcal Y^{n,\alpha+1}_t.
	\end{equation}
\end{enumerate}
In both cases, the first mode should be KPZ, while the second should be diffusive.

We present the study of fluctuations for the case I, but the strategy and the technical tools developed therein are the same for both the subcase and for the general case. 
Briefly, the argument consists in: proving tightness of the sequences $\{\mathcal Z_t^n\}_{n\in\N}$, $\{\widetilde{\mathcal Z}^n_t\}_{n\in\N}$, with respect to the Skorohod topology of ${D}([0,T],{\mathcal D'(\mathbb T)})$; characterizing the limit points of these sequences: we have to prove that any limit point of $\{\mathcal Z_t^n\}_{n\in\N}$ and of $\{\widetilde{\mathcal Z}^n_t\}_{n\in\N}$ is an energy solution --- see Section~\ref{sec:SES}, to the SBE and the OUE respectively.

For the sake of completeness, we present the multispecies version of the definition of energy solution~\cite{BFS,GO}. We say that the process $(\mathcal{Y}^\alpha_\cdot,\mathcal{Y}^{\alpha+1}_\cdot)\in{D}([0,T],{\mathcal D}'(\mathbb T)^2)$ is a \emph{multi-species energy solution} of a system of coupled SBEs if:

\begin{itemize}
	\item[\bf{i)}] For any $t\in[0,T]$, $(\mathcal{Y}^\alpha_t,\mathcal{Y}^{\alpha+1}_t)$ is a spatial white noise with covariance 
	\begin{equation}
	{\rm{Cov}}(\partial_u\eta^\alpha(f),\partial_u\eta^{\alpha+1}(g))=\Gamma_{x,y}^{\alpha,\alpha+1}\int f(u)g(u)du,
	\end{equation}
	where $\Gamma^{\alpha,\alpha+1}_{x,y}$ denotes the covariance
\begin{equation}
\Gamma^{\alpha,\alpha+1}_{x,y}:=\int_{\Omega_n}\bar\xi^\alpha_x\bar\xi^{\alpha+1}_y\nu(d\eta)
\end{equation}
and $\Gamma^\alpha_x=\Gamma^{\alpha,\alpha}_{x,x}$ the variance;
\item[\bf{ii)}] The process $(\mathcal{Y}^\alpha_\cdot,\mathcal{Y}^{\alpha+1}_\cdot)$ satisfy a $\mathbb L^2$-energy condition: there exist a process\\ $(\mathcal{A}^\alpha_\cdot, \mathcal{A}^{\alpha+1}_\cdot) \in C([0,T],{\mathcal D'(\mathbb{T})^2})$ whose action on $f,g\in {\mathcal D(\mathbb{T})}$ is the uniformly $\mathbb L^2$-Cauchy limit $(\mathcal{A}^\alpha_\cdot(f),\mathcal{A}^{\alpha+1}_\cdot(g))$. We also impose that $\mathcal{A}^{\alpha(+1)}_\cdot(f)$ has zero quadratic variation, $\langle\mathcal{A}^{\alpha(+1)}(f)\rangle_t=0$;
	\item[\bf{iii)}] For any $f\in {\mathcal D(\mathbb{T})}$, $$\mathcal{M}^{\alpha(+1)}_t(f)=\mathcal{Y}^{\alpha(+1)}_t(f)-\mathcal{Y}^{\alpha(+1)}_0(f)-\int \mathcal{Y}^{\alpha(+1)}_s(\Delta f)ds-\mathcal{A}^{\alpha(+1)}_t(f)$$ is a Brownian motion of variance $\Gamma^{\alpha}_x t\|\nabla f\|^2_{2}$ with respect to the filtration generated by $\mathcal{Y}^{\alpha(+1)}_\cdot$;
	\item[\bf{iv)}] The time-reversed processes $\hat{\mathcal{Y}}^{\alpha(+1)}_\cdot=\mathcal{Y}^{\alpha(+1)}_{T-\cdot}$ and $\hat{\mathcal{A}}^{\alpha(+1)}_\cdot=\mathcal{Y}^{\alpha(+1)}_T-\mathcal{Y}^{\alpha(+1)}_{T-\cdot}$ also satisfy \textbf{iii)}.
\end{itemize}

To prove tightness we will rely on widely used criteria: the Mitoma's criterion~\cite{mitoma}, that reduces tightness of distribution-valued processes to tightness of real-valued processes; the Aldous' criterion; and the Kolmogorov--Prohorov--Centsov's criterion, invoked in the case of processes with continuous paths.

\subsection{{Convergence to the SBE}}\label{sec:convSBE}
According to \eqref{eq:modes2}, first we look at the field $\mathcal Z_t^n=\mathcal Y^{n,\alpha}_t$. Recall \eqref{eq:densityfieldNL} and, so far, take there  $v_\alpha(n)=0$, this means there is no velocity for this mode. We will see below what is the velocity to take in such a way that we can control the terms in the martingale expansion. 
We denote by $\bar\xi^\alpha_x=\xi^\alpha_x-{\mathbbm{E}}_\nu[\xi^\alpha_x]$ the centred variable.
The infinitesimal current (centred by adding and subtracting appropriately $\frac 13\xi^{\alpha}_{x(+1)}$) is
\begin{equation}\label{eq:current1}
j^\alpha_{x,x+1}=\bar{\xi}^\alpha_{x+1}-\bar{\xi}^\alpha_x-\frac{E}{n^\gamma}\bar{\xi}^\alpha_x\bar{\xi}^\alpha_{x+1}+\frac{E}{6n^\gamma}(\bar{\xi}^\alpha_x+\bar{\xi}^\alpha_{x+1}).
\end{equation}
Note that for the choice of the rates in the case I, the current above does not depend on the other species of the system. 
We recall the definition of discrete derivatives~\eqref{eq:discretederivatives}. Dynkin's formula for $\mathcal Z_t^n$ reads as
\begin{equation}\label{eq:DynkinZ3}
\begin{aligned}
&\mathcal M_t^n(f)={\mathcal Z}^{n}_t(f)-{\mathcal Z}^{n}_0(f)-\int_0^t ds \mathcal Z_s^n(\Delta_n  f)+\frac{\sqrt n}{n^{\gamma}}E\int_0^t ds \sum_{x\in\mathbb{T}_n}\nabla^+_n f\left(\tfrac{x}{n}\right)\bar\xi^\alpha_{x}(sn^2)\bar\xi^\alpha_{x+1}(sn^2)\\
&-\frac{\sqrt n}{n^{\gamma}}\frac{E}{3}\int_0^t ds \sum_{x\in\mathbb{T}_N}\nabla^+_n f\left(\tfrac{x}{n}\right)\bar\xi^\alpha_{x}(sn^2)+\frac{1}{n^{\gamma+1/2}}\frac{E}{3}\int_0^t ds \sum_{x\in\mathbb{T}_n}\Delta^+_n f\left(\tfrac{x}{n}\right)\bar\xi^\alpha_{x}(sn^2).
\end{aligned}
\end{equation}
{For $\gamma>0$, the variance of the last term  goes to zero with $n\to\infty$.}  We  also get rid of the remaining degree one term on the last line by  redefining the field $\mathcal{Y}^{n,\alpha}_t$ with a moving frame of velocity $v_\alpha(n)=-\frac E 3n^{2-\gamma}$.
To conclude, we need a second order Boltzmann--Gibbs principle to {write $\bar\xi^\alpha_{x}\bar\xi^\alpha_{x+1}$ in terms of a quadratic functional of the density field, in the case $\gamma=1/2$; and to show that the term with the product  $\bar\xi^\alpha_{x}\bar\xi^\alpha_{x+1}$ vanishes as $n\to+\infty$, when $\gamma>1/2$. Therefore, from now on take $\gamma=1/2$.} If we denote by $\overrightarrow\xi^{\epsilon n,\alpha}_{x}$ the average of $\bar{\xi}^\alpha_x$ on a box of size $\epsilon n$ to the right of $x$ as in \eqref{eq:boxes_averages}, the Boltzmann--Gibbs principle allows to replace the degree-two term with $(\overrightarrow\xi^{\epsilon n,\alpha}_{x}(sn^2))^2$ and therefore, we can close the equation in terms of the field ${\mathcal Z}^{n}$ as
\begin{equation}\label{eq:martingaleZ3}
\begin{aligned}
\mathcal M_t^n(f)={\mathcal Z}^{n}_t(f)-{\mathcal Z}^{n}_0(f)-\int_0^t ds  \mathcal Z_s^n(\Delta_nf)
-E\int_0^t ds{\frac 1n} \sum_{x\in\mathbb{T}_n}\nabla^+_n T^+_{v_\alpha(n)s}f\left(\tfrac{x}{n}\right)\left(\mathcal{Z}^n_s({i}_\epsilon(\tfrac{x+v_\alpha (n)s}{n}))\right)^2,
\end{aligned}
\end{equation}
with $ i_\epsilon(x)=\frac{1}{\epsilon}\Id_{[x,x+\epsilon]}$ and $v_\alpha(n)=-\frac E 3n^{3/2}.$
{By making a change of variables we rewrite the last expression as 
\begin{equation}\label{eq:martingale_new}
\begin{aligned}
\mathcal M_t^n(f)={\mathcal Z}^{n}_t(f)-{\mathcal Z}^{n}_0(f)-\int_0^t ds  \mathcal Z_s^n(\Delta_nf)-E\int_0^t ds{\frac 1n} \sum_{x\in\mathbb{T}_n}\nabla^+_n f\left(\tfrac{x}{n}\right)\left(\mathcal{Z}^n_s({i}_\epsilon(\tfrac{x}{n}))\right)^2.
\end{aligned}
\end{equation}
}

\subsection{{Convergence to OUE}}\label{sec:convOU}
We now discuss the second mode $\widetilde{\mathcal Z}_t^n=\mathcal Y^{n,\alpha}_t+2\mathcal Y^{n,\alpha+1}_t$. For the species  $\alpha+1$ the infinitesimal current is
\begin{equation}\label{eq:current2}
j^{\alpha+1}_{x,x+1}=\bar{\xi}^{\alpha+1}_{x+1}-\bar{\xi}^{\alpha+1}_x-\frac{E}{n^\gamma}(\bar{\xi}^\alpha_x\bar{\xi}^{\alpha+1}_{x+1}+\bar{\xi}^{\alpha+1}_{x}\bar{\xi}^\alpha_{x+1})-\frac{E}{6n^\gamma}(\bar{\xi}^\alpha_x+\bar{\xi}^{\alpha}_{x+1}+\bar{\xi}^{\alpha+1}_{x}+\bar{\xi}^{\alpha+1}_{x+1}).
\end{equation}
We take the linear combination of the currents~\eqref{eq:current1} and~\eqref{eq:current2}, and we redefine the fields $\mathcal Y^{N,\alpha}_t$, $\mathcal Y^{N,\alpha+1}_t$ by choosing moving frames of velocity {$\tilde v_\alpha(n)=\tilde v_{\alpha+1}(n)=\tfrac E3 n^{a-\gamma}$}. {We observe that in this case both fields have the same velocity, nevertheless the field  $\mathcal Y^{N,\alpha}_t$ is taken now in a different frame than the one we analysed in the previous subsection. }  In this way, the degree one terms in the current disappear and Dynkin's formula gives 
\begin{equation}\label{eq:DynkinZ4}
\begin{aligned}
\widetilde{\mathcal M}_t^n(f)=&\widetilde{\mathcal Z}^{n}_t(f)-\widetilde{\mathcal Z}^{n}_0(f)-\int_0^t ds \widetilde{\mathcal Z}_s^n(\Delta_nf)\\&+{E\frac{\sqrt n}{n^\gamma}}\int_0^t ds \sum_{x\in\mathbb{T}_n}\nabla^+_nT^+_{\tilde v_{\alpha+1}( n)s}f\left(\tfrac{x}{n}\right)[\bar{\xi}^{\alpha}_x\bar{\xi}^{\alpha}_{x+1}+\bar{\xi}^{\alpha}_x\bar{\xi}^{\alpha+1}_{x+1}+\bar{\xi}^{\alpha+1}_x\bar{\xi}^{\alpha}_{x+1}](sn^2).
\end{aligned}
\end{equation}
{From the second order Boltzmann-Gibbs principle, the last term in last display vanishes as $n\to+\infty$ if $\gamma>1/2$. But for $\gamma=1/2$, and again from that principle,  we are able } to replace $\bar\xi^{\alpha(+1)}_x$ with its average on a box of size $\epsilon n$ to the left of $x$, $\overleftarrow\xi^{\epsilon n,\alpha(+1)}_{x}$, and $\bar\xi^{\alpha(+1)}_{x+1}$ with its average on a box of size $\epsilon n$ to the right of $x+1$, $\overrightarrow\xi^{\epsilon n,\alpha(+1)}_{x}$, so that we can close the equation in terms of the fields $\widetilde{\mathcal Z}^{n}$ and ${\mathcal Z}^{n}$ as
\begin{equation}\label{eq:DynkinZ4BG}
\begin{aligned}
\widetilde{\mathcal M}_t^n(f)=\widetilde{\mathcal Z}^{n}_t(f)-\widetilde{\mathcal Z}^{n}_0(f)-\int_0^t ds \widetilde{\mathcal Z}_s^n(\Delta_nf)
+E\int_0^t ds {\frac 1n}\sum_{z\in\mathbb{T}_n}\nabla^+_n f\left(\tfrac{z}{n}\right){\mathcal Z}^n_s ( i_\epsilon(\tfrac{z+(\tilde v_\alpha(n)-v_{\alpha}(n))s}{n}))\widetilde{\mathcal Z}^n_s( i_\epsilon(\tfrac zn)).
\end{aligned}
\end{equation}

Here we performed the change of variables $x=z-{\tilde v}n^{3/2}s$.   This means  that to  conclude the argument  we have to prove that} the last integral of~\eqref{eq:DynkinZ4BG} goes to zero with $n\to\infty$. Informally, this should follow from the fact that we are evaluating the field $\mathcal{Z}^n_s$ in a frame which is not its ``natural'' one: taking the $n\to\infty$ limit, the support of $\Ie$ goes to $-\infty$, away from the support of the test function $f$. A rigorous proof of this fact will appear in~\cite{GO}.

\thanks

\end{document}